\documentclass[3p,preprint,12pt]{elsarticle}

\usepackage{lineno}
\usepackage{amsmath}
\usepackage{amsfonts}
\usepackage{amssymb}
\usepackage{isomath}
\usepackage{enumerate}
\usepackage{booktabs}
\usepackage{multirow}
\usepackage{pifont}
\usepackage{enumitem}
\usepackage{algpseudocode}
\usepackage[english]{babel}
\usepackage{amsmath}
\usepackage[dvips]{epsfig}
\usepackage{graphicx,amsmath,amsthm,amsfonts}
\usepackage{epsfig}
\usepackage{amsopn}
\usepackage{color}
\usepackage{amsmath,amssymb,amsfonts,stmaryrd}
\usepackage{subfig}
\usepackage{graphicx}
\usepackage{tikz, pgfplots}
\usepackage{hhline}
\usepackage{float}
\usepackage{multirow}
\usepackage{appendix}
\usepackage{placeins}
\usepackage{xparse}
\usepackage{mathtools}
\mathtoolsset{showonlyrefs=true}
\usepackage{fixltx2e}
\MakeRobust{\eqref}
\usepackage{everypage}
\usepackage{xspace}
\usepackage{ifthen}
\usepackage{enumitem}
\usepackage{url}
\reversemarginpar
\usepackage{algorithm,algpseudocode}
\definecolor{notes}{RGB}{0,255,178}
\usepackage{array}
\usepackage{amsmath}
\usepackage{amssymb}
\usepackage{amsbsy}
\usepackage{graphicx}
\usepackage{float}
\usepackage{xcolor}
\usepackage{bbm}
\usepackage{gensymb}
\usepackage{graphicx}
\usepackage{subfig}
\usepackage{multirow}
\usepackage{enumitem}
\usepackage{caption}
\DeclareMathOperator*{\argmax}{arg\,max}
\DeclareMathOperator*{\argmin}{arg\,min}

\NewDocumentCommand{\normDG}{O{\cdot} O{j} O{}}{\ensuremath{\|\ifthenelse{\equal{#1}{}}{\cdot}{#1}\|_{DG,\ifthenelse{\equal{#2}{}}{j}{#2}}\ifthenelse{\equal{#3}{}}{}{^{#3}}}\xspace}
\NewDocumentCommand{\normL}{O{\cdot} O{2} O{\Om} O{}}{\ensuremath{\|\ifthenelse{\equal{#1}{}}{\cdot}{#1}\|_{L^{\ifthenelse{\equal{#2}{}}{1}{#2}}(\ifthenelse{\equal{#3}{}}{\Om}{#3})}\ifthenelse{\equal{#4}{}}{}{^{#4}}}\xspace}
\NewDocumentCommand{\normH}{O{\cdot} O{1} O{\Om} O{}}{\ensuremath{\|\ifthenelse{\equal{#1}{}}{\cdot}{#1}\|_{H^{\ifthenelse{\equal{#2}{}}{1}{#2}}(\ifthenelse{\equal{#3}{}}{\Om}{#3})}\ifthenelse{\equal{#4}{}}{}{^{#4}}}\xspace}
\NewDocumentCommand{\norms}{O{\cdot} O{s} O{j} O{}}{\ensuremath{\ltrivert\ifthenelse{\equal{#1}{}}{\cdot}{#1}\|_{\ifthenelse{\equal{#2}{}}{s}{#2},\ifthenelse{\equal{#3}{}}{j}{#3}}\ifthenelse{\equal{#4}{}}{}{^{#4}}}\xspace}
\NewDocumentCommand{\Aa}{O{j} O{\cdot} O{\cdot}}{\ensuremath{\mathcal{A}_{\ifthenelse{\equal{#1}{}}{j}{#1}}(\ifthenelse{\equal{#2}{}}{\cdot}{#2},\ifthenelse{\equal{#3}{}}{\cdot}{#3})}\xspace}
\NewDocumentCommand{\mcal}{O{} O{} O{}}{\ensuremath{\mathcal{#1}\ifthenelse{\equal{#2}{}}{}{_{#2}}\ifthenelse{\equal{#3}{}}{}{^{#3}}}\xspace}

\newcommand{\Om}{\ensuremath{\Omega}\xspace}

\newcommand{\xx}{\boldsymbol}
\newcommand{\yy}{\displaystyle}

\journal{Computers in Biology and Medicine}

\begin{document}

\begin{frontmatter}

\title{Integration of activation maps of epicardial veins in computational cardiac electrophysiology}
\author[label1]{Simone Stella}
\ead{simone.stella@polimi.it}

\author[label2]{Christian Vergara\corref{cor1}}
\ead{christian.vergara@polimi.it}
\cortext[cor1]{Corresponding author}

\author[label3]{Massimiliano Maines}
\ead{massimiliano.maines@apss.tn.it}

\author[label3]{Domenico Catanzariti}
\ead{domenico.catanzariti@apss.tn.it}

\author[label1]{Pasquale Claudio Africa}
\ead{pasqualeclaudio.africa@polimi.it}

\author[label3]{Cristina Demattè}
\ead{cristina.dematte@apss.tn.it}

\author[label4]{Maurizio Centonze}
\ead{maurizio.centonze@apss.tn.it}

\author[label5]{Fabio Nobile}
\ead{fabio.nobile@epfl.ch}

\author[label3]{Maurizio Del Greco}
\ead{Maurizio.DelGreco@apss.tn.it}

\author[label6]{Alfio Quarteroni}
\ead{alfio.quarteroni@polimi.it}

\address[label1]{MOX, Dipartimento di Matematica, Politecnico di Milano, Milan, Italy}
\address[label2]{LABS, Dipartimento di Chimica, Materiali e Ingegneria Chimica "Giulio Natta", Politecnico di Milano, 
Milan, Italy} 
\address[label3]{Divisione di Cardiologia, Ospedale S. Maria del Carmine, Rovereto (TN), Italy}
\address[label4]{U.O. di Radiologia di Borgo-Pergine, Ospedale di Borgo Valsugana, Borgo Valsugana (TN), Italy}
\address[label5]{CSQI, Institute of Mathematics, \'Ecole Polytechnique Fédérale de Lausanne, Switzerland}
\address[label6]{MOX, Dipartimento di Matematica, Politecnico di Milano, Milan, Italy \& Institute of Mathematics, \'Ecole Polytechnique Fédérale de Lausanne, Switzerland (professor emeritus)}

%%%%%%%%%%%%%%                              ABSTRACT                         %%%%%%%%%%%%%%%%%%%%%%%%
\begin{abstract} 
	\noindent In this work we address the issue of validating the monodomain equation
	used in combination with the Bueno-Orovio ionic model for the prediction of the activation
	times in cardiac electro-physiology of the left ventricle. To this aim, we consider four patients who suffered from Left 
	Bundle Branch Block (LBBB). We use activation maps performed at the septum as input data for the model
	and maps at the epicardial veins for the validation. In particular, a first set (half) of the latter
	are used to estimate the conductivities of the patient and a second set (the remaining half) to compute
	the errors of the numerical simulations. We find an excellent agreement between measures and numerical 
	results. Our validated
	computational tool could be used to accurately predict activation times at the epicardial veins
	with a short mapping, i.e. by using only a part (the most proximal) of the standard 
	acquisition points,
	thus reducing the invasive procedure and exposure to radiation. 
	\end{abstract}
%%%%%%%%%%%%%%%%%%%%%%%%%%%%%%%%%%%%%%%%%%%%%%%%%%%%%%%%%%%%%%%%%%%%%%%%%%%%%%%%%%%%%%%%%%%%%%%%%%%%%


\end{frontmatter}

%%%%%%%%%%%%%%                            INTRODUCTION                       %%%%%%%%%%%%%%%%%%%%%%%%
\section{Introduction}
Electrophysiology is a fundamental research field in applied mathematics since Hodgkin 
and Huxley described for the first time the propagation of action potentials in cells 
\cite{Hodgkin_1952}. Mathematical and numerical modeling in cardiac electrophysiology in the last 
decades have assumed a key role to better understand cardiac muscle function and to study how several 
cardiac diseases develop and to provide concrete answers to clinical problems. 
Challenging issues consist of selecting
accurate and efficient numerical methods for the approximate solution of such models and 
in the estimation of the model parameters to fit patient-specific data.

To describe the propagation of the electrical signal in the heart muscle two possible mathematical approaches are 
the bidomain and the monodomain models \cite{Colli_Franzone, Colli_Franzone_2005}. The former is the result of the 
application of conservation of charge together with constitutive models. It describes the propagation 
of the trans-membrane potential and compute both the internal and external cell potentials. The 
monodomain model is a simplification of the bidomain one, assuming that the 
external and internal conductivity tensors are proportional. Both models have 
to be coupled with a ionic model, a system of ODEs which describes the evolution of the 
trans-membrane potential in a single cell by means of suitable gating variables.

Validation of these mathematical models, i.e. the certification that a model could be 
used for predictive purposes, is fundamental to provide clinicians a reliable tool to study 
and predict the cardiac function accurately. 

A first attempt in this direction consisted of comparing numerical results
(obtained either by the bidomain or monodomain model) with measures of electrical activity obtained by 
optical imaging during \textit{ex-vivo} experiments on animal hearts. In particular, qualitative comparisons have been
provided in \cite{Potse_2006,Krishnamoorthi_2014} for action potentials and
in \cite{Bordas_2011,Deng_2015} for activation times. Other works instead
quantified the discrepancy between results and measures; in particular in \cite{Muzikant_2002,Barone_2020} the authors 
focused on the action potential, whereas in \cite{Camara_2011,Relan_2011,Wang_2013} on activation times. 
In \cite{Pop_2013}, the authors considered a similar calibration of the monodomain model
by using data of activation times acquired by catheters in pigs.

A second set of results focused on human ideal geometries, proposing benchmark simulation protocols to be validated against gold-standard activation times obtained from experimental measures. In particular,  
\cite{Niederer_2011} considered a slab of tissue, whereas \cite{Niederer_2010} an idealized
left ventricle.

Another group of studies addressed qualitative comparisons in real human geometries.
In particular, in \cite{Ten_Tusscher_2009} the authors compared
the phase distribution during ventricular fibrillation with some reference electrical data,
whereas in \cite{Augustin_2016} they studied the reliability of ECG obtained with an electro-mechanical simulation.
Instead, in \cite{Konukoglu_2011}, the authors performed a calibration of the monodomain model on human data
by using a probabilistic model (Bayesian influence method).

Other works have considered the inclusion of electrical data of human activation time obtained by catheter 
in a computational framework to calibrate the conduction properties of the tissue \cite{Chinchapatnam_2008,Lines_2009,Dossel_2012,Sermesant_2012,vergarap2,Palamara_2014,Sanchez_2017,Lee_2019}.
For example, in \cite{vergarap2,Palamara_2014} a quantitative comparison between in-vivo measures of activation
times acquired at the endocardium and those provided by the Eikonal equation in the presence of a personalized Purkinje network was carried out. Instead, in \cite{Corrado_2018} the authors validated the monodomain model with
a Mitchell-Schaeffer ionic model in the left atrium under the S1 and S2 stimulation protocols.   
Among the previous studies, only \cite{Lee_2019} considered measures at the ventricular epicardial veins.      
In particular, the reaction-Eikonal equation was used to find the most accurate model among six electrophysiology surrogate models able to reproduce the electrical activation during right ventricle apex stimulation.

Of previous cited works, the majority conducted a {\sl calibration} of the model considered
for their study, i.e. all the data at disposal have been used to estimate some conductivity parameters to match the data.
Instead, only few works used different sets of data (in some cases 
also of different nature) to calibrate the model and to assess the error against clinical measures 
(with a little abuse of terminology, we refer to this as {\sl cross-validation}, as done in statistics when 
the results obtained from the statistical analysis are compared with the testing data selected randomly and not
used in the calibration).
In particular, the two sets are given by action potential data in pigs
at different stimulation cycle lengths in \cite{Barone_2020}, 
depolarization times and $dp/dt$ together with the blood pressure 
in \cite{Sermesant_2012}, activation times at different locations in \cite{vergarap2,Palamara_2014} for the Eikonal model, ECG and epicardial activation times in \cite{Lee_2019} for the reaction-Eikonal model, atrial activation times for different stimulations in \cite{Corrado_2018}.

Our work aims at proposing an approach for validating the monodomain model 
in the context of the electrical propagation in the 
human left ventricle using patient-specific activation measures acquired at the epicardial veins. 
In particular, we only consider excitations with sinus rhythm, and not, tachicardia, fibrillation,
and stimulated cases. Notice that, with sinus rhythm here we refer also to pathological cases
	(such as Left Bundle Branch Block (LBBB) without scar), provided that the front propagates without irregular or chaotic patterns. 
For each patient we had at disposal activation times at some points located
at the septum and at the epicardial veins. 
The first dataset was used to provide a patient-specific input for the numerical simulation, whereas the second one 
was split into two subgroups, the first one used to estimate the conductivities of the patient
and the second one to validate the accuracy of the numerical solution.
Notice that the protocol of validation we described is quite general in the sense 
that the use of measurements acquired by means of catheters is not 
essential; in principle, comparison and minimization against 
activation maps acquired for example by electrode array recordings and optical mapping are possible.
At the best of authors knowledge, this is the first cross-validation test for the monodomain model performed 
for the ventricular activation against clinical measures of activation times.

%%%%%%%%%%%%%%%%%%%%%%%%%%%%%%%%%%%%%%%%%%%%%%%%%%%%%%%%%%%%%%%%%%%%%%%%%%%%%%%%%%%%%%%%%%%%%%%%%%%%%

%%%%%%%%%%%%%%                             SECTION 2                         %%%%%%%%%%%%%%%%%%%%%%%%
\section{Methods}
\subsection{Mathematical and numerical model}
In this section we present the mathematical model for cardiac electrophysiology, \textit{i.e.} the monodomain model,
and the corresponding numerical methods considered in this work.
The choice of the monodomain model has been motivated by its lower computational cost,
yet providing comparable accuracy, with respect to the bidomain
model, at least in the cases of sinus rhythms \cite{Potse_2006}.

\begin{figure}[H]
	\centering
	\includegraphics[scale=0.18]{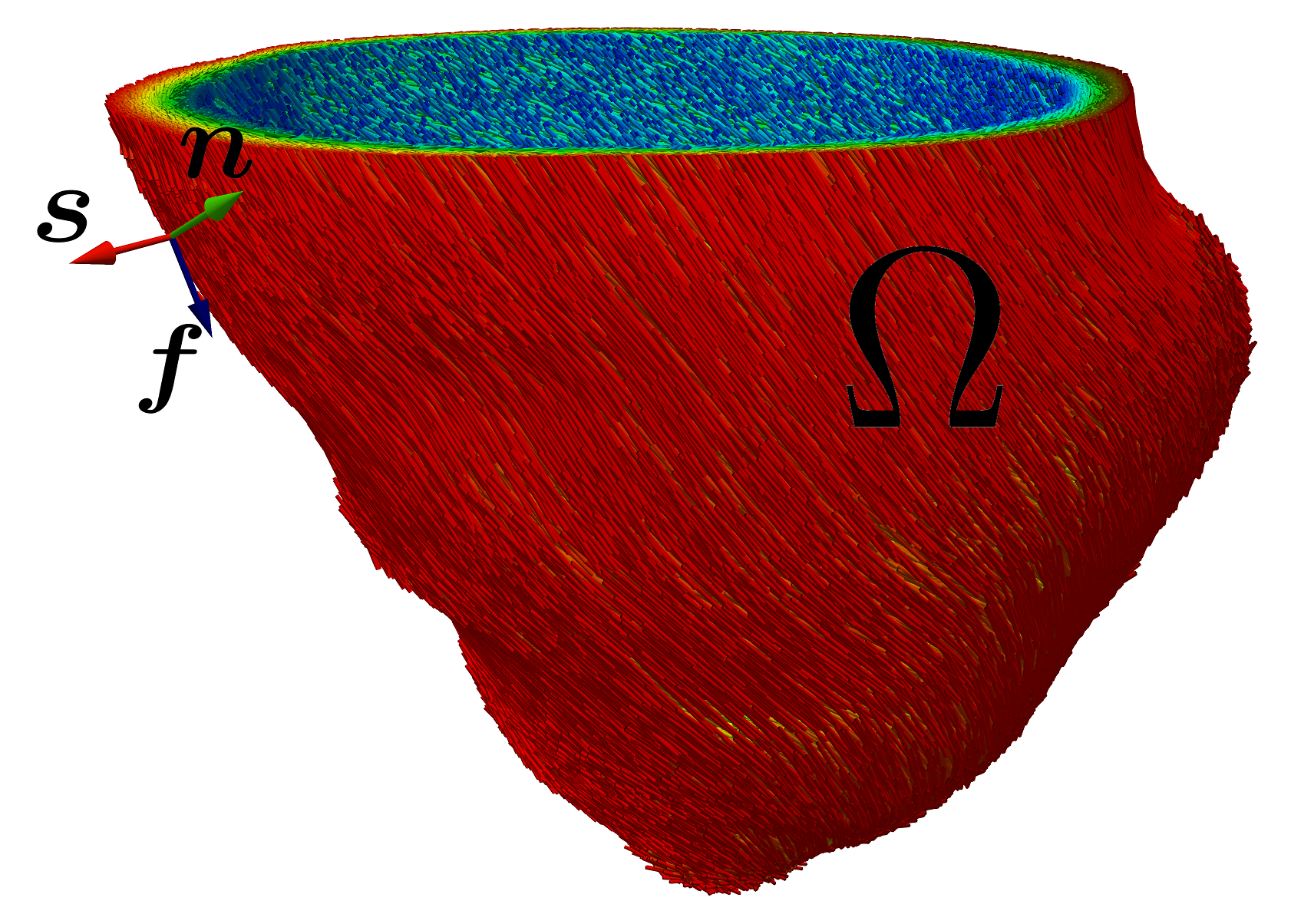}
	\caption{Computational domain.} \label{fig:domain}
\end{figure}

Referring to Figure \ref{fig:domain}, the cardiac tissue conductivity is modeled as a tensor $\mathbb{D}$ defined as
\begin{equation}\label{eq:conduct_tensor}
	\begin{alignedat}{1} \mathbb{D} = \sigma_s \mathbbm{1} + (\sigma_f - \sigma_s)\boldsymbol f \otimes 
		\boldsymbol f + (\sigma_n - \sigma_s)\boldsymbol n \otimes \boldsymbol n, \end{alignedat}
\end{equation}
where $\sigma_f$, $\sigma_n$ and $\sigma_s$ are the conductivities along the fibers direction $\xx f$, the normal direction $\xx n$
and the transversal  direction $\xx s$ (orthogonal to the sheets plane $\xx f$-$\xx n$), respectively, 
to be determined in order to fit the patient-specific activation times.

Assuming an external applied stimulus $I_\mathrm{app}$ (provided in our case by activation time measures
located at the septum, see below), the monodomain model in the computational domain \(\Omega\) and over a time interval \((0, T]\) reads \cite{Colli_Franzone_2014,bookIheart}:

Find, for each $t$, the electrical potential $u:\Omega \to\mathbb R$ and the gating variables $\boldsymbol w:\Omega \to\mathbb{R}^3$, such that
\begin{subequations}
	\label{eq:monodomain_system}
	\begin{align}
		& \chi C_m\frac{\partial u}{\partial t} - \nabla\cdot\left(\mathbb{D}
		\nabla u\right) + \chi I_\mathrm{ion}(u, \boldsymbol w) = I_\mathrm{app} & 
		\text{in } \Omega \label{eq:monodomain},\\ 
		%& \mathbb{D}\nabla u\cdot\mathbf{\nu} = 0, & \text{on } \partial \Omega \times (0, T],\\
		& \displaystyle\frac{d\boldsymbol w}{dt} = \boldsymbol R(u,\boldsymbol w) \label{eq:ionicmodel} & \text{in } \Omega,
	\end{align}
\end{subequations}
where the ionic current $I_\mathrm{ion}(u,\boldsymbol w)$ and $\boldsymbol R\in\mathbb{R}^3$ are chosen according to the
Bueno-Orovio ionic model \cite{Bueno-Orovio2008}. In particular, $I_\mathrm{ion}=I_\mathrm{ion,1}+I_\mathrm{ion,2}+I_\mathrm{ion,3}$, with

\begin{subequations}
	\label{eq:BO_ionic_model}
	\begin{align}
		&I_\mathrm{ion,1} = \yy- \frac{H(u-V_1)(u- V_1)(\widetilde V -u)w_1}{\tau_1},\\[2.5ex]
		&I_\mathrm{ion,2} = \yy\frac{1-H(u-V_2)(u-V_o)}{H(u-V_o)(\tau_{o2}-\tau_{o1})+\tau_{o1}} +
		\frac{H(u-V_2)}{H(u-V_2)(\tau_{22}-\tau_{21})+\tau_{21}},\\[2.5ex]
		&I_\mathrm{ion,3} = -\yy\frac{H(u-V_2)}{\tau_3}w_2w_3,
	\end{align}
\end{subequations}
for suitable constants $V_o,\, V_1,\, V_2,\, \widetilde V,\, \tau_1,\, \tau_3,\, \tau_{o1},\, \tau_{o2},\, \tau_{21},\, \tau_{22}$ and where $H$ is the Heaviside function.
Moreover, $\chi$ is the surface area-to-volume ratio and $C_m$ the trans-membrane capacitance. 
System \eqref{eq:monodomain_system} has been equipped with suitable initial conditions for $u$ and $\boldsymbol w$ and homogeneous Neumann conditions on $\partial \Omega$ for $u$.

The time discretization we used to numerically approximate system \eqref{eq:monodomain_system} relies
at each time step $t^{n+1}=(n+1)\Delta t$ first on a forward Euler method for 
equation \eqref{eq:ionicmodel}, \textit{i.e.} 
\begin{equation}
	\label{eq:ionicmodel_dt}
	\yy\frac{\xx w^{n+1}-\xx w^n}{\Delta t} = \xx R(u^n,\xx w^n)\qquad \text{in } \Omega,
\end{equation}
\(\Delta t\) being the time step length, and then on a first order semi-implicit method for \eqref{eq:monodomain}, i.e.
\begin{equation}
	\label{eq:monodomain_dt}
	\chi C_\mathrm{m}\frac{u^{n+1}-u^n}{\Delta t} - \nabla\cdot\left(\mathbb{D}
	\nabla u^{n+1}\right) + \chi \left(I^{n+1}_\mathrm{ion,1}+I^{n+1}_\mathrm{ion,2}+I^{n+1}_\mathrm{ion,3}\right)
	= I_\mathrm{app}(t^{n+1}) 
	\qquad\text{in } \Omega,
\end{equation}
with
\begin{equation}
	\label{eq:monodomain_dt2}
	\left.
	\begin{array}{l}
		I^{n+1}_\mathrm{ion,1} = \yy- \frac{H(u^n-V_1)(u^{n+1}-V_1)(\widetilde V -u^n)w^{n+1}_1}{\tau_1},\\[2.5ex]
		I^{n+1}_\mathrm{ion,2} = \yy\frac{1-H(u^n-V_2)(u^{n+1}-V_o)}{H(u^n-V_o)(\tau_{o2}-\tau_{o1})+\tau_{o1}} +
		\frac{H(u^n-V_2)}{H(u^n-V_2)(\tau_{22}-\tau_{21})+\tau_{21}},\\[2.5ex]
		I^{n+1}_\mathrm{ion,3} = -\yy\frac{H(u^n-V_2)}{\tau_3}w^{n+1}_2w^{n+1}_3,
	\end{array}
	\right.
\end{equation}
where the diffusion term has been treated implicitly and the non-linear terms have been linearized,
see also \cite{pathmanathan2012computational}. In particular, at each time we first updated $\boldsymbol{w}^{n+1}$ for a given $u^n$ by means of \eqref{eq:ionicmodel_dt}, then we solved 
\eqref{eq:monodomain_dt} and \eqref{eq:monodomain_dt2} by using the up-to-date gating variables.
Such time discretization lead to a conditionally stable method with a bound on time step $\Delta t$ which is independent of the mesh size and which is milder with respect to the $\Delta t$ required to reach the desired accuracy \cite{Fernandez_2010,Colli_Franzone_2014}.

As for the space discretization, we used continuous Finite Elements of order 1 (Q1) on hexahedral meshes. 
The ionic current term \(I_\mathrm{ion}\) has been discretized using the Ionic Current Interpolation (ICI) method \cite{krishnamoorthi2013numerical}: first \(I_\mathrm{ion}\) has been computed using the values of \(u\) and \(\boldsymbol{w}\) at the degrees of freedom, then it has been interpolated at quadrature nodes. Such approach is relatively inexpensive and less memory-demanding than solving the ODE system and computing \(I_\mathrm{ion}\) directly at quadrature nodes (SVI), while the numerical accuracy is not affected at the small mesh size required to capture the propagating front.
The resulting linear system arising at each time step has been solved by the 
GMRES method \cite{saad1986gmres} preconditioned with the Jacobi preconditioner.

We have used a time step $\Delta t = 0.025ms$ 
and a characteristic mesh size $h\simeq 0.35\,mm$. Notice that both these values are small enough to recover an accurate propagation front \cite{Fernandez_2010,Gurev_2010,Colli_Franzone_2014,Arevalo_2016},
see in particular \cite{Hurtado_2017} for the case of hexahedral meshes with Q1 Finite Elements,
as in our case.

The values of parameters used in \eqref{eq:monodomain_dt}-\eqref{eq:monodomain_dt2} are reported in Table \ref{table:parameters} \cite{Bueno-Orovio2008}.
\begin{table}[H]
	\centering
	\subfloat{
		%\footnotesize
		\centering
		\begin{tabular}{|c|c|c|c|c|c|c|}
			\hline
			$C_m$ $[F/m^2]$ & $\chi$ $[m^{-1}]$ & $V_o$ $[s]$ & $V_1$ $[s]$ & $V_2$ $[s]$ & $\widetilde V$ $[s]$ & $\tau_1$ $[s]$\\
			\hline
			$0.01$ & $1\times 10^5$ & $0.006$ & $0.3$ & $0.015$ & $1.58$ & $11\times 10^{-3}$\\
			\hline
		\end{tabular}
	}\\
	\subfloat{
		%\footnotesize
		\centering
		\begin{tabular}{|c|c|c|c|c|}
			\hline
			$\tau_3$ $[s]$ & $\tau_{o1}$ $[s]$ & $\tau_{o2}$ $[s]$ & $\tau_{21}$ $[s]$ & $\tau_{22}$ $[s]$ \\
			\hline
			$2.8723\times 10^{-3}$ & $6\times 10^{-3}$ & $6\times 10^{-3}$ & $43\times 10^{-3}$ & $0.2\times 10^{-3}$\\
			\hline
		\end{tabular}
	}\caption{Values of the coefficients used in \eqref{eq:monodomain_dt}-\eqref{eq:monodomain_dt2}.}\label{table:parameters}
\end{table}

Finally, we highlight the importance of including cardiac fibers in electrophysiology models,
since electrical propagation occurs in a different way along the fibers and orthogonally to them. 
Since standard imaging techniques do not provide geometric information on the fibers,
whose dimension is typically smaller than the spatial resolution, the fiber orientation over the myocardial tissue has been here determined 
using the Laplace-Dirichlet rule-based algorithm described in \cite{Bayer_2012}. 
	In particular, we used a linear rule with the following boundary values for the fibers and the sheets angles: $-60^{\circ}$ 
	for fibers at the epicardium, $60^{\circ}$ for fibers at endocardium, $20^{\circ}$ for sheets on epicardium and $-20^{\circ}$ for sheets on endocardium \cite{Bayer_2012,Doste_2019}.

The Laplace-Dirichlet rule-based fiber generation algorithm as well as all the numerical methods for the monodomain equation have been implemented within \texttt{life\textsuperscript{x}} \\
(https://lifex.gitlab.io/lifex), a new in-house developed high-performance \texttt{C++} library mainly focused on cardiac applications, based on the \texttt{deal.II} Finite Element core \cite{dealII91}; 
for further details
on the implementation of the fibers generation see \cite{Piersanti_2020}. 
Electrophysiology simulations ran on 48 cores of a 192 cores node Platinum Intel® Xeon® 8160 @2.1GHz with 1.7TB RAM.
%%%%%%%%%%%%%%%%%%%%%%%%%%%%%%%%%%%%%%%%%%%%%%%%%%%%%%%%%%%%%%%%%%%%%%%%%%%%%%%%%%%%%%%%%%%%%%%%%%%%%

%%%%%%%%%%%%%%                                                    %%%%%%%%%%%%%%%%%%%%%%%%
\subsection{Processing of geometric and electrical data}\label{sec:processing}
In this section we illustrate the strategies used for data processing.
In particular, in Section \ref{sec:geom} we detail the geometric reconstruction, whereas in Section 
\ref{sec:map} the activation time maps and their integration with geometric data.
Finally, in Section \ref{sec:invpb} we set the inverse problem used to provide the integration
of electrical data into the numerical experiments and the estimation of the conductivities 
for a validation of the monodomain model in the context of a sinus rhythm. All the medical data, both geometric and electrical,
have been provided by Ospedale S. Maria del Carmine, Rovereto (TN), Italy.

Four patients have been considered, from now on referred to as P1, P2, P3, and P4.
They were all affected by a LBBB, a cardiac conduction abnormality due to an 
interruption of the 
electrical conduction in the His bundle which causes a delayed activation of the left ventricle. LBBB was diagnosed owing to ECG evaluation.
In all patients scar regions were absent and the mapping procedure was performed by clinicians 
	during the implantation of the Cardiac Resynchronization Therapy (CRT), in accordance with the standard clinical procedures.

All patients have been previously informed and gave their full consent for the acquisition of both the geometric and electrical data and for the successive mathematical analyses.

%%%%%%%%%%%%%%%%%%%%%%%%%%%%%%%%%%%%%%%%%%%%%%%%%%%%%%%%%%%%%%%%

\subsubsection{Imaging data acquisition and geometric reconstruction}\label{sec:geom}
The four subjects considered in the present study underwent a cineMRI steady-state free precession ECG/Retro with a 1.5-Tesla MRI Unit (Magnetom Aera, Siemens Medical Systems, Erlangen, Germany). The following parameters have been used: in-plane resolution $1.7 \times 1.3\, mm^2$ and slice thickness $8 mm$; 
TR (repetition time) $38.64$ $ms$; TE (echo time) $1.41$ $ms$; flip angle $79\degree$; averages $1$; SNR (signal noise ratio) $1$.

\begin{figure}[H]
	\centering
	\subfloat[\textbf{P1}]{\includegraphics[width=0.45\textwidth]{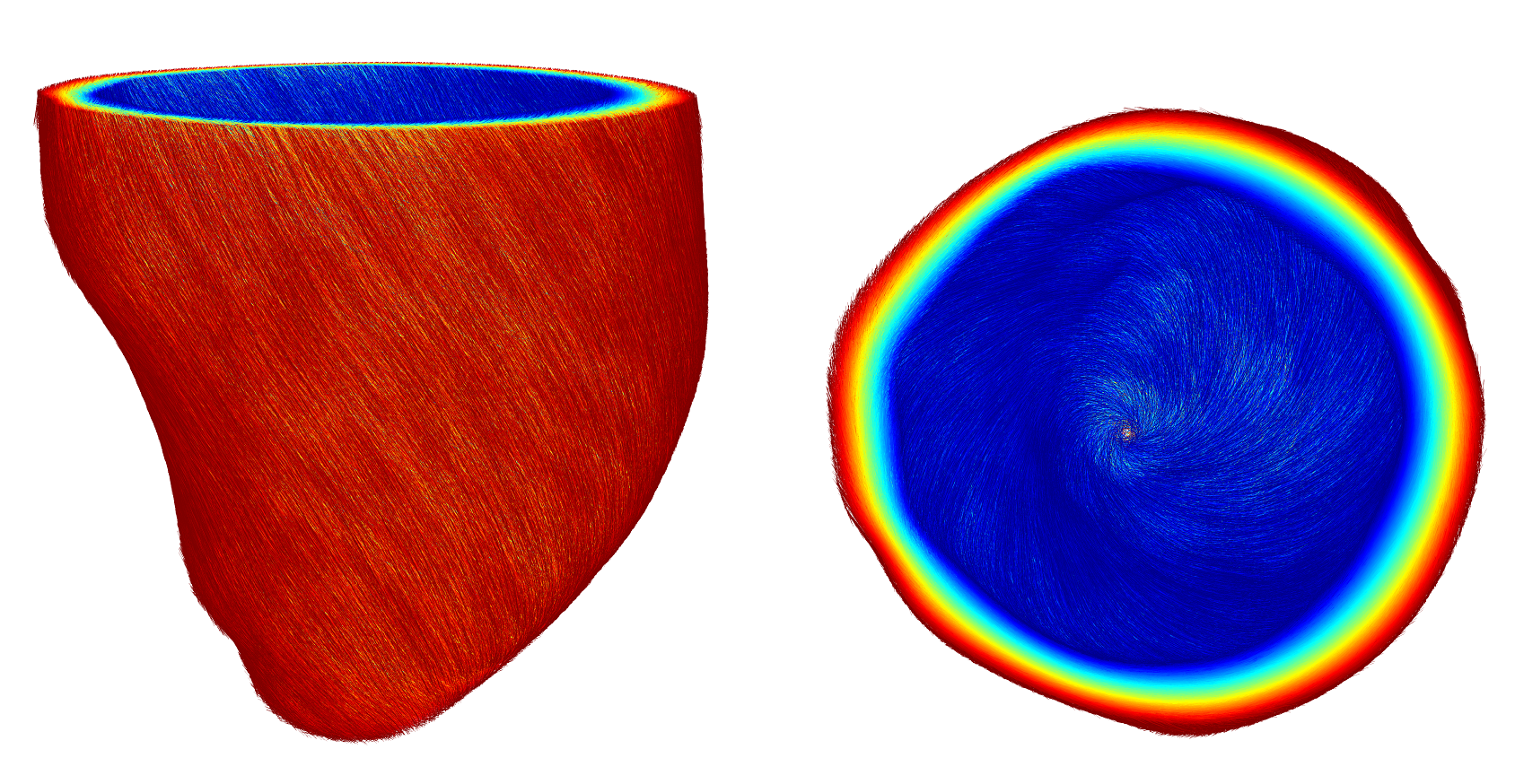}}\hfill
	\subfloat[\textbf{P2}]{\includegraphics[width=0.45\textwidth]{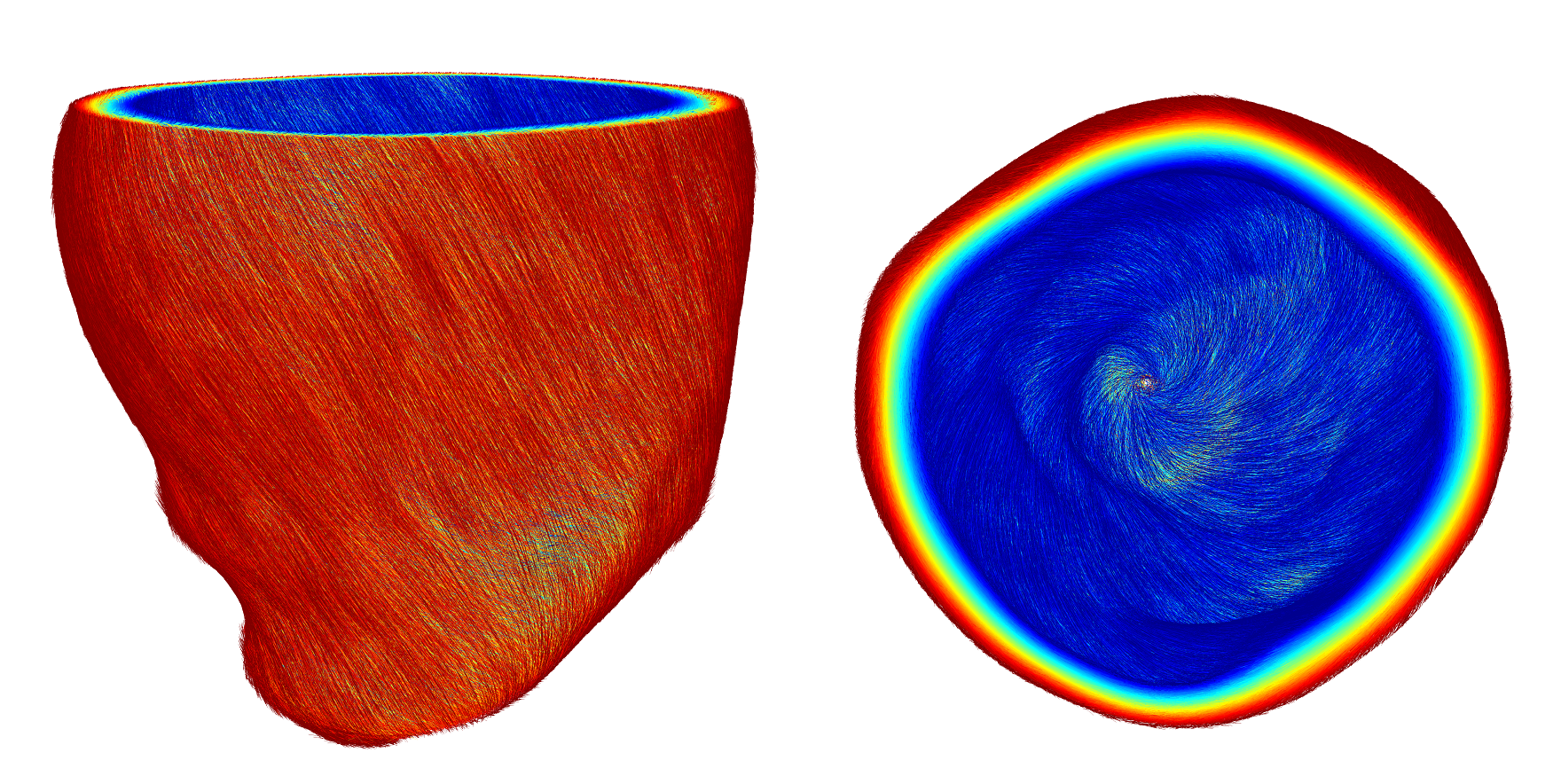}}\\
	\subfloat[\textbf{P3}]{\includegraphics[width=0.45\textwidth]{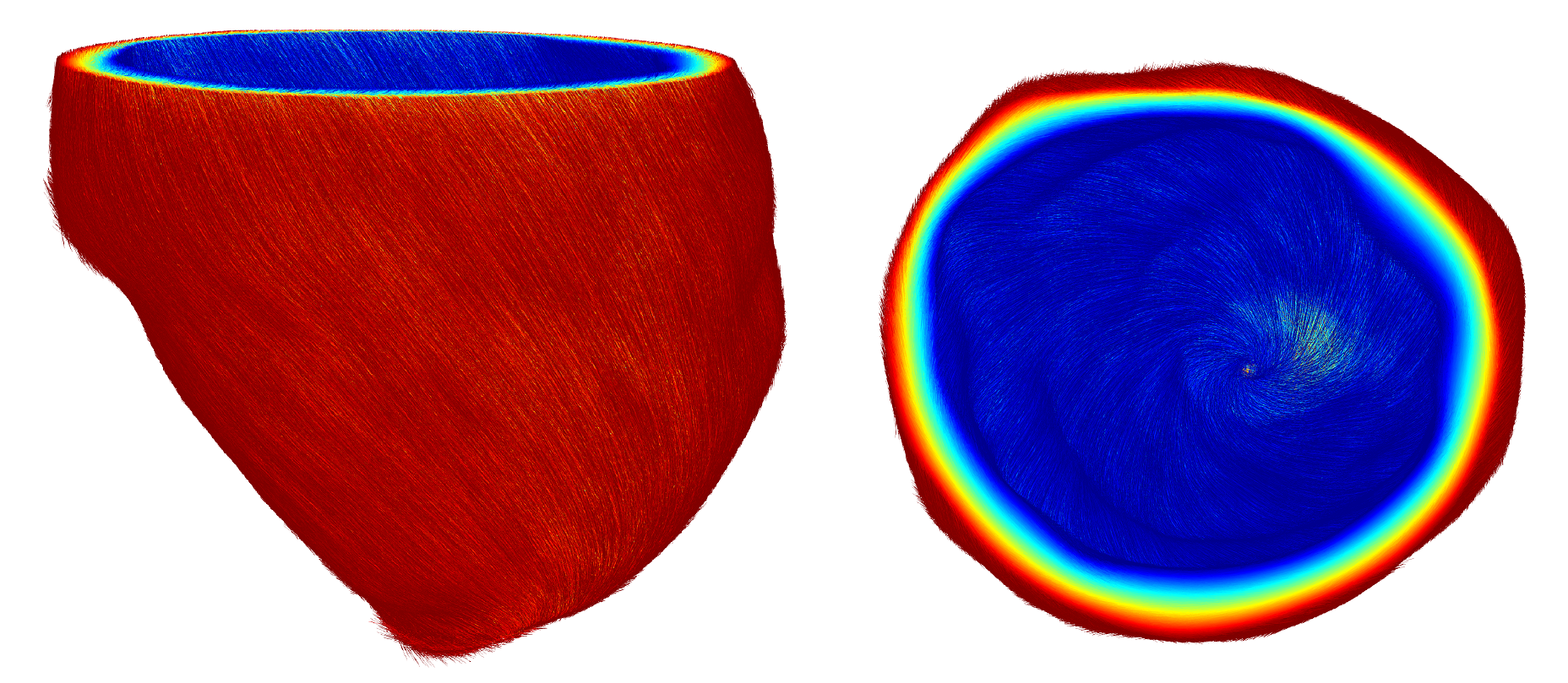}}\hfill
	\subfloat[\textbf{P4}]{\includegraphics[width=0.45\textwidth]{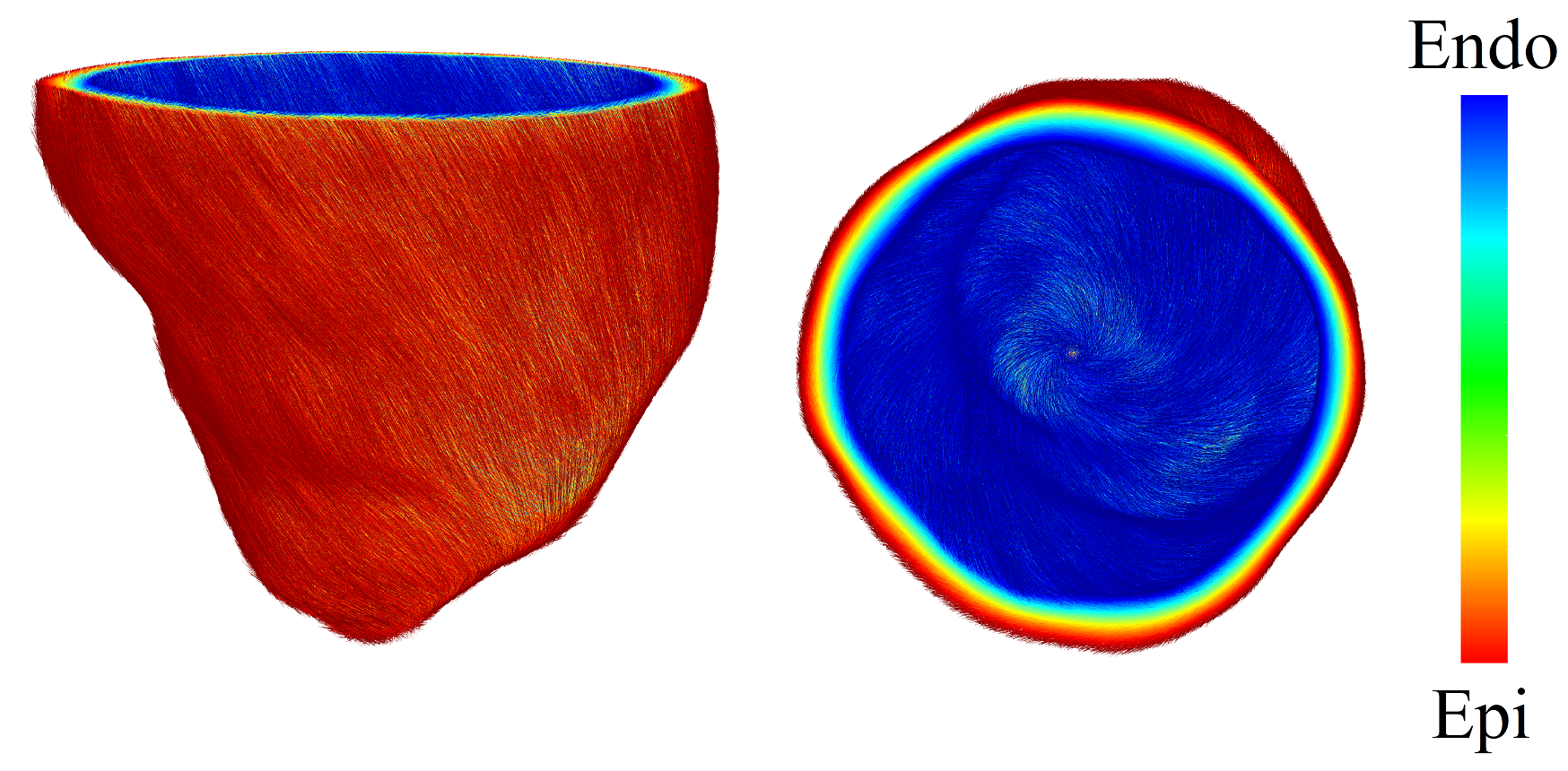}}
	\caption{Patient-specific reconstructed geometries. For each case: on the left the frontal view, on the right base view. Muscle fibers generated by Laplace-Dirichlet rule-based algorithm are integrated into the 3D model.} \label{fig:geom}
\end{figure}

We performed a semi-automatic segmentation of the left ventricle geometry by using the free open-source software  \texttt{MITK} \cite{Mitk}, which allowed us to segment each MRI slice and to interpolate the ventricle surface between slices. In particular, for each patient MRI dataset, we have preliminary applied a manual correction of slice misalignment due to motion and breath artefacts.
The next step of the segmentation procedure consisted of capping the two surfaces at the base of the ventricle, then to connect them with a triangulated base using the approach reported in \cite{fedele2019}.
Once two closed triangulated surfaces of both epicardium and endocardium of the left ventricle were obtained, a surface mesh was generated by means of a set of new meshing tools \cite{fedele2019} developed as an extension to the \texttt{VMTK} software \cite{Antiga_2008}. A remeshing procedure on the whole closed surface was then performed in order to prescribe a target mesh size. Finally,
a volumetric hexahedral mesh of the left ventricle was generated. In Figure \ref{fig:geom} the four patient-specific reconstructed volumes of the left ventricles with the generated muscle fibers are displayed.

\captionsetup[subfigure]{labelformat=empty}

\subsubsection{Activation maps acquisition}\label{sec:map}
A mapping of the septal surface located in the right ventricle and of the coronary veins located at the epicardium of the left ventricle has been performed (from 7 to 10 days after the MRI acquisition) for all the four subjects through the use of the  {\sl Ensite Precision} system \cite{Eitel_2010} to record local activation times
(referred to also as electrical data).
The latter represent, for each point of acquisition, the time difference between two instants, one measured on an extracardiac reference electrode and another one representing the steepest negative intrinsic deflection in the electrogram recorded on an intracardiac electrode and indicating that the activation wavefront has reached the point under investigation. The Ensite Precision system allows to perform an accurate real-time three dimensional catheter navigation to obtain maps of activation times. For this study, a 5 Fr steerable 10-pole catheter has been inserted through the left subclavian vein. 
These data are available from multiple beats in order to filter out anomalous beats
such as extrasystoles and take signals with similar morphology.
In particular, the system records signals if their duration is in the range $(900-1100)*HR/60\,ms$,
with $HR$ the heart rate of the patient, and discard those ones that feature a morphology which is not aligned
with the others. 
In this way, the inter-beat variability in our activation recordings is limited, thus yielding a signal as much as possible homogeneous.

We subdivided the electrical data at disposal into three subsets,
each of them composed by the activation times and the coordinates of the corresponding points:
\begin{itemize}
	\item[-] {\sl Septal data}: acquired at the septum and
	used here to provide the input current $I_\mathrm{app}$ in \eqref{eq:monodomain_dt};
	\item[-] {\sl epicardial veins data, group I}: corresponds to half of the epicardial vein 
	measures we have at disposal, in particular it is composed by the points with the earliest activation times.
	This set is used to calibrate the conductivities of the patient; 
	\item[-] {\sl epicardial veins data, group II}: corresponds to the remaining half of epicardial vein 
	measures, in particular it is composed by the points with the highest activation times.
	In the spirit of a cross-validation, this set has been used to compute the discrepancies with the numerical solution and thus provide a validation of the latter. 
\end{itemize}

The choice of using the earliest and latest activated points for the calibration and validation groups,
respectively, allowed us to maximize the ''separation'' between the two groups and thus to strengthen the
validation which has been performed on a distinct set with respect to the one used for calibration. This could also have important clinical consequences as highlighted in the Discussion.  

In Table \ref{table:measurements} we report the number of activation time measurements $N^S$ (septal),
$N^V_{I}$ (epicardial vein, group I) and 
$N^V_{II}$ (epicardial vein, group II) used in this work for each patient.
\begin{table}[H]
	\centering
	\begin{tabular}{c|c|c|c|c|}
		\cline{2-5}
		& P1 & P2 & P3 & P4 \\
		\hline
		\multicolumn{1}{ |c|  }{$N^S$} & 38 & 15 & 9 & 4 \\
		\hline
		\multicolumn{1}{ |c|  }{$N^V_I$} & 8 & 26 & 19 & 16 \\
		\hline
		\multicolumn{1}{ |c|  }{$N^V_{II}$} & 8 & 26 & 18 & 15 \\
		\hline
	\end{tabular}
	\caption{Total number of measurements acquired at the septum and at the epicardial veins for each patient.}\label{table:measurements}
\end{table}

In order to include in our simulation framework measured activation maps obtained from Ensite Precision
(electrical data) onto the reconstructed geometries obtained from MRI (geometric data), we needed to merge geometric and electrical data. Since the MRI units and the Ensite Precision are two distinct systems collecting clinical data, the reconstructed patient-specific geometry and the corresponding activation map point cloud were linked to two distinct reference systems.

In order to make them compatible, we applied to each patient the following procedure based on the following three steps:
\begin{itemize}
	\item[-]  \textbf{Reference points selection}: we selected three points for each set of data (geometric and electrical) as a reference.
	Two of them have been chosen on the coronary sinus, the third one on the septal surface of the right ventricle (see Figure \ref{fig:2}, left block). We then verified, owing to the clinicians experience, that the two points
	of each couple in fact corresponded to the same physical point;
	\item[-] \textbf{Geometric alignment}: we applied a rotation and translation to the point cloud of electrical data 
	so that the three couples of reference points identified at the previous step coincided (see Figure \ref{fig:2}, middle block). 
	In fact, this guaranteed that the point cloud of electrical data lay as much as possible 
	in correspondence of the geometric MRI data;
	\item[-] \textbf{Nearest Neighbor Search (NNS) projection}: 
	as a consequence of the geometric discrepancy still present between the two sets of data (geometric and electrical ones) due to the different instants 
	acquisitions of MRI and Ensite Precision, we applied for the point cloud 
	of electrical data the Euclidean NNS procedure \cite{Andrews_2001}. This allowed us to select for each point of the cloud the 
	nearest one among the geometric data belonging to the external (epicardial or septal) surface of the left ventricle. We finally moved each cloud point accordingly, so that it coincided with the geometric nearest one determined at the previous step (see Figure \ref{fig:2}, right block).
\end{itemize}
This procedure, depicted in Figure \ref{fig:2}, has been implemented in a Python script interacting with the Paraview software (Reference points selection and geometric alignment) and in MATLAB (NNS projection).

\begingroup
\centering
\begin{figure}[H]
	\centering
	\includegraphics[scale=0.5]{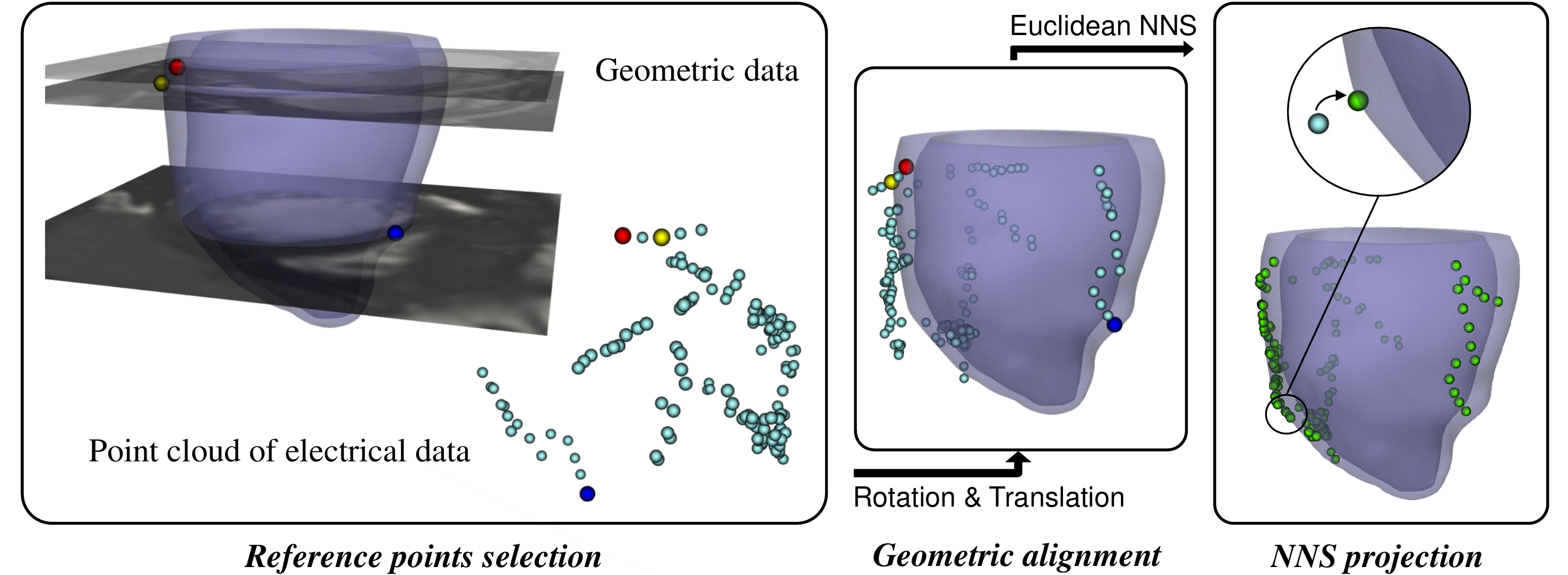}
	\caption{Graphical representation of the procedure for alignment of electrical and geometric data.} \label{fig:2}
\end{figure}
\endgroup

At the end of the alignment procedure we could identify 
\begin{itemize}
	\item[-] the septal data 
	$(\xx x^S_1,\tau^S_1),$
	$\ldots,$
	$(\xx x^S_j,\tau^S_j),$
	$\ldots,$
	$(\xx x^S_{N^S},\tau^S_{N^S})$;  
	\item[-] the epicardial veins data belonging to group I $(\xx x^{V-I}_1,\tau^{V-I}_1),$
	$\ldots,$
	$(\xx x^{V-I}_j,\tau^{V-I}_j),$
	$\ldots,$
	$(\xx x^{V-I}_{N^V_I},\tau^{V-I}_{N^V_I})$;
	\item[-] the epicardial veins data belonging to group II $(\xx x^{V-II}_1,\tau^{V-II}_1),$
	$\ldots,$
	$(\xx x^{V-II}_j,\tau^{V-II}_j),$
	$\ldots,$
	$(\xx x^{V-II}_{N^V_{II}},\tau^{V-II}_{N^V_{II}})$.
\end{itemize}

In Table \ref{table:AT} we report some information on the location of the points, 
in particular the minimum and maximum activation time for each of the groups (I and II)
and the distance between them, intended as the minimum geodesic distance among all the possible couples
of points belonging to the two groups. 
\begin{table}[H]
\centering
\begin{tabular}{|c||c|c||c|c||c|}
\hline
&&&&&\\
Patient & min $AT^{I}$ [$s$] & max $AT^{I}$ [$s$] & min $AT^{II}$ [$s$] & max $AT^{II}$ [$s$] &  $D^{I-II}$ [$cm$] \\
\hline
P1 & 0.110 & 0.152 & 0.157 & 0.179 & 0.17\\
\hline
P2 & 0.068 & 0.105 & 0.106 & 0.122 & 0.11\\
\hline
P3 & 0.107 & 0.119 & 0.121 & 0.142 & 0.09\\
\hline
P4 & 0.077 & 0.106 & 0.110 & 0.148 & 0.12\\
\hline
\end{tabular}\caption{Values of the minimum and maximum activation time measurements for group I and group II, and   geodesic distance $D^{I-II}$ between the two groups.}\label{table:AT}
\end{table}

\begin{figure}[H]
	\centering
	\subfloat[\textbf{P1}]{\includegraphics[width=0.45\textwidth]{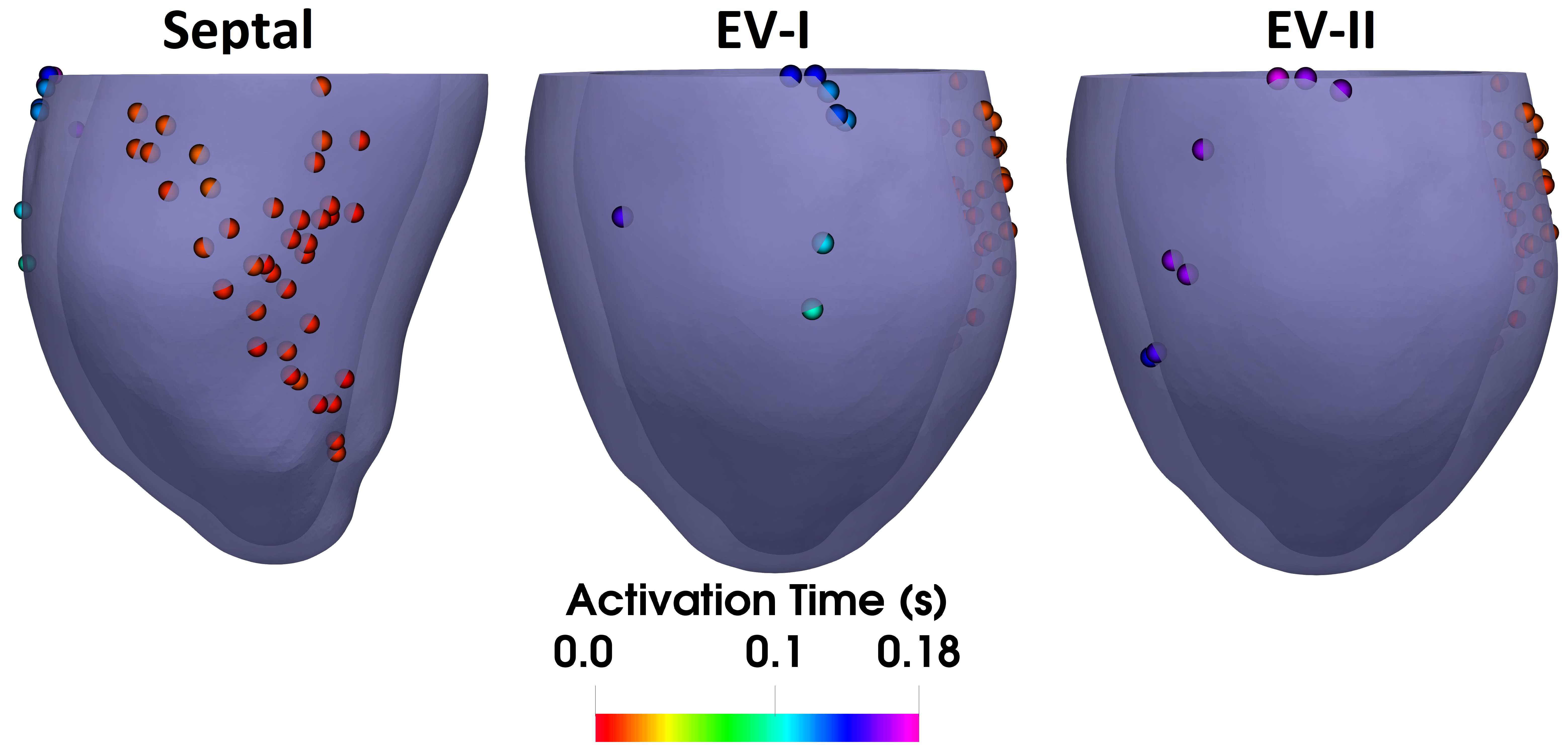}}\hfill
	\subfloat[\textbf{P2}]{\includegraphics[width=0.45\textwidth]{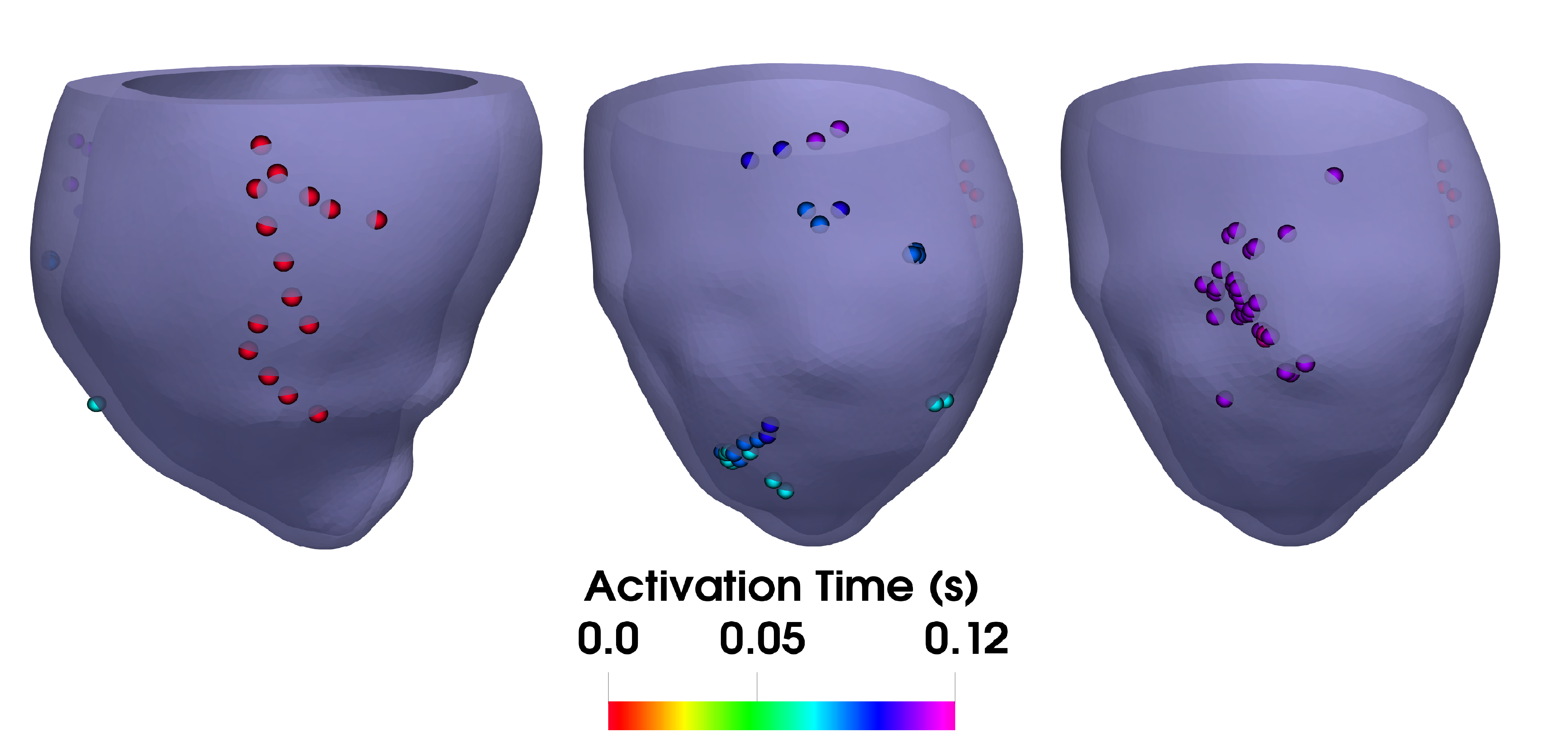}}\\
	\subfloat[\textbf{P3}]{\includegraphics[width=0.45\textwidth]{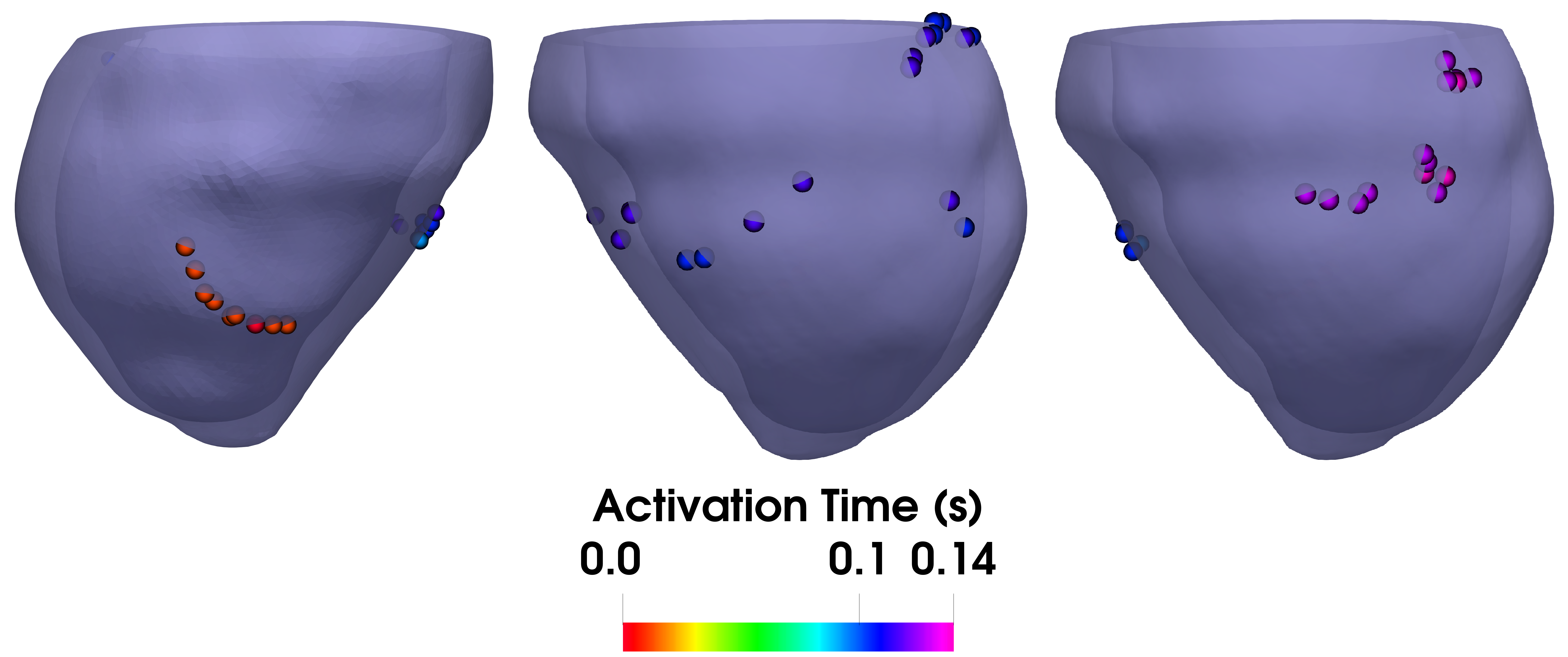}}\hfill
	\subfloat[\textbf{P4}]{\includegraphics[width=0.45\textwidth]{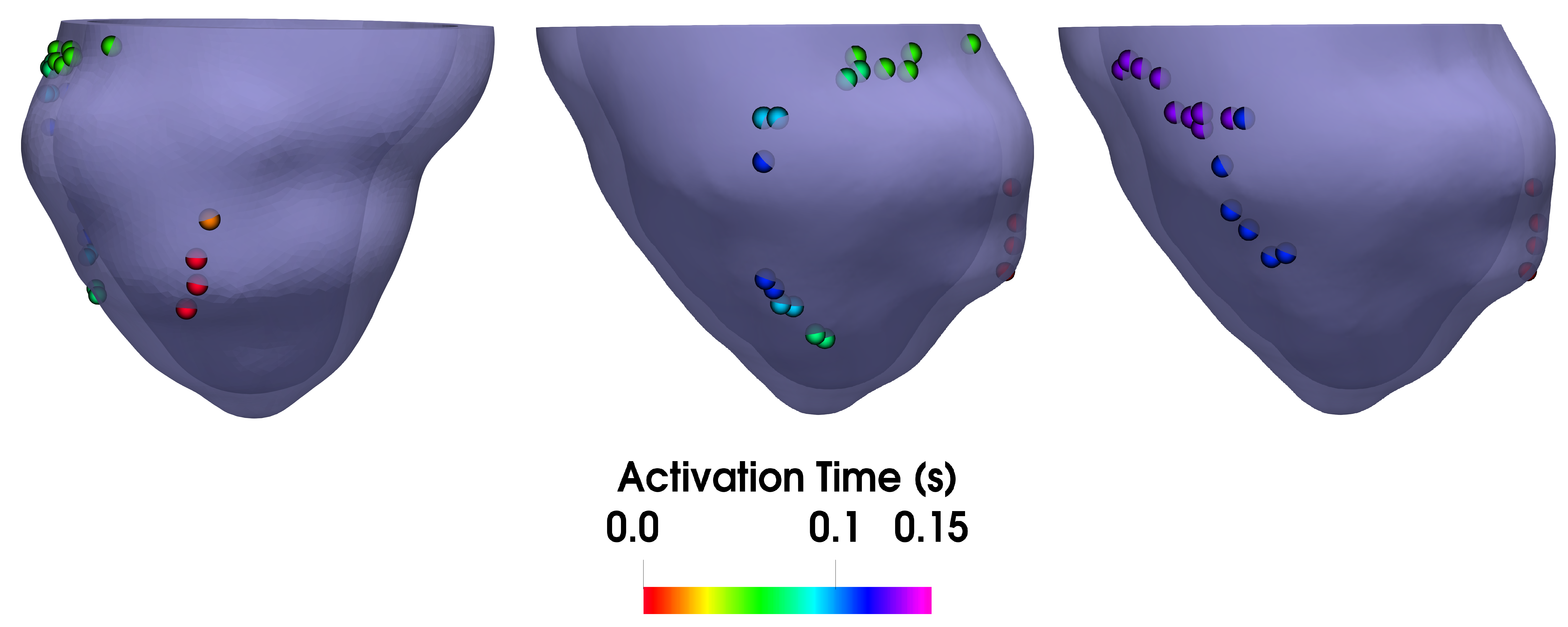}}
	\caption{Maps of activation time after the alignment with the geometric MRI data.
		For each case: Left, septal data; Middle, epicardial veins data - group I (EV-I); 
		Right, epicardial veins data - group II (EV-II).} \label{fig:3}
\end{figure}

In Figure \ref{fig:3}, we report the final results of the alignment procedure for all the four patients with the 
corresponding patient-specific activation time measurements.

\captionsetup[subfigure]{labelformat=empty}

\subsubsection{Estimation of conductivities}\label{sec:invpb}  % Validation through minimization problem

We decided to use the 
septal measurements to define the applied current $I_\mathrm{app}$ in the Finite Elements approximation of problem 
\eqref{eq:monodomain_dt}, whereas the epicardial veins measurements were used 
to estimate the conductivities $\sigma_f,\,\sigma_s$ and $\sigma_n$ in \eqref{eq:conduct_tensor}
and to assess a significant step towards the validation in the context of sinus rhythm
of the monodomain model equipped with the Bueno-Orovio ionic model in a patient-specific context. 
Since in this work we considered patients without any scar, we assumed without loss
of generality that the three conductivites are constant in space and time 
\cite{Sebastian_2008, Bayer_2012, Corrado_2015, Potse_2006}. 

The applied current was built accordingly as follows:
\begin{equation}
	I_\mathrm{app}(\xx x,t) = \begin{cases} 112500\mbox{ }\mu A\mbox{ }cm^{-3} & \mbox{if } 
		(\xx x,t)=(\xx x_j^S,\tau_j^S)\quad \mbox{for some} \,\, j=1,\ldots,N^S,
		\\ 0\mbox{ }\mu A\mbox{ }cm^{-3} & \mbox{elsewhere.}\end{cases}
\end{equation}
The value of the applied current was chosen as the lowest value able to 
allow the electrical signal to propagate in the ventricle \cite{Niederer_2011}.

In what follows, we specified how the activation times \(\tau^{h-I}_j\) and \(\tau^{h-II}_j\)
at the two groups of epicardial veins points were computed from the numerical simulations.
In particular, 
\(\tau^{h-I}_j\) (resp. \(\tau^{h-II}_j\)) has been defined as the discrete time instant where the Finite Elements approximation of the trans-membrane potential $u_h$\footnote{With a slight abuse of notation, we denote by $u_h$ the Finite Elements time discretized solution $u_h^n$.} at the computational point \(\xx x_j^{V-I}\) (resp. \(\xx x_j^{V-II}\)) varies at its highest rate, \textit{i.e.}
\begin{equation}
	\label{eq:tau}
	\tau^{h-\beta}_j = t^{\bar{n}}, \quad \text{where } \bar{n}= \argmax_{n} \left|\frac{u_h^n\left(\xx x_j^{V-\beta}\right) - u_h^{n-1}\left(\xx x_j^{V-\beta}\right)}{\Delta t}\right|,\quad\beta=I,II,
\end{equation}
where a first order Euler approximation has been used, consistently with the order of the 
monodomain time discretization.

In order to maximize the agreement between numerical simulations and clinical measurements, we looked for 
the conductivities $\xx\sigma=(\sigma_f,\sigma_s,\sigma_n)$ in the physiological range 
$\Sigma=$
$(0.70,2.20)$
$\times$
$(0.16,0.48)$
$\times$
$(0.03,0.10))$
$k\Omega^{-1}cm^{-1}$  \cite{Clerc_1976,Roberts_1979,Roberts_1982,Le_Guyader_2001,Stinstra_2005}.
Specifically, we wanted to minimize the discrepancy between the computed activation times
$\tau^{h-I}_j$ and the epicardial veins measures belonging to group I, $\tau^{V-I}_j$.
To this aim, we introduced the following discrete functional:
\begin{equation}\label{eq:functional}
	F(u_h(\xx\sigma))=\sum_{j=1}^{N^{V-I}}\frac{1}{2}\left|\tau^{h-I}_j(u_h(\xx\sigma))-\tau^{V-I}_j\right|^2.
\end{equation}
Notice that we have highlighted the 
dependence of $u_h$ on $\xx\sigma$. 

The optimization problem then reads: Find the optimal value $\widehat{\xx\sigma}$ such that 
\begin{equation}\label{eq:minimization}
	\widehat{\xx\sigma}=\argmin_{\xx\sigma\in\Sigma} F(u_h(\xx\sigma)),
\end{equation}
subjected to the Finite Elements approximation of the discretized-in-time monodomain problem \eqref{eq:ionicmodel_dt}-\eqref{eq:monodomain_dt}-\eqref{eq:monodomain_dt2}.

To solve minimization problems similar to the previous one, 
some efficient strategies have been proposed for example in \cite{Yang_2017} for synthetic data 
and \cite{Barone_2020} for optical measurements on animal hearts. In \cite{Costa_2013} a simpler iterative method has been proposed on a slab geometry, exploiting the proportional relation between conduction velocites and conductivities.
Here, since we considered a minimization problem with \textit{in-vivo} human electrical data, 
we 
	preferred to use all patient-specific electrical data (group I) at disposal in order to have a more
	robust result. We followed a basic {\sl direct search} method
which is robust with respect to the noise of the electrical data. This is based on starting by 
an initial guess of $\xx\sigma$ taken in $\Sigma$ and on ongoing corrections obtained by
solving the monodomain problem and by evaluating the functional \eqref{eq:functional}. 
In particular, given suitable acceleration parameters $\beta_f=0.45,\,\beta_s=0.1,\,\beta_n=0.05$, the mean error
$E^{(k)}=\sum_{j=1}^{N^{V-I}}\left(\tau^{h-I}_j(u^{(k)}_h(\xx\sigma))-\tau^{V-I}_j\right)$ at iteration $k$ is computed.
The new value of each component of $\xx\sigma$ is then updated as follows:
\[
\sigma^{(k+1)}_\gamma = \sigma^{(k)}_\gamma + \beta_\gamma E^{(k)}   ,\quad \gamma=f,s,n.
\]  

%%%%%%%%%%%%%%%%%%%%%%%%%%%%%%%%%%%%%%%%%%%%%%%%%%%%%%%%%%%%%%%%%%%%%%%%%%%%%%%%%%%%%%%%%%%%%%%%%%%%%

%%%%%%%%%%%%%%                             SECTION 4                         %%%%%%%%%%%%%%%%%%%%%%%%
\section{Results}

\subsection{Standard scenario}
\begin{figure}[H]
	\centering
	\includegraphics[scale=0.3]{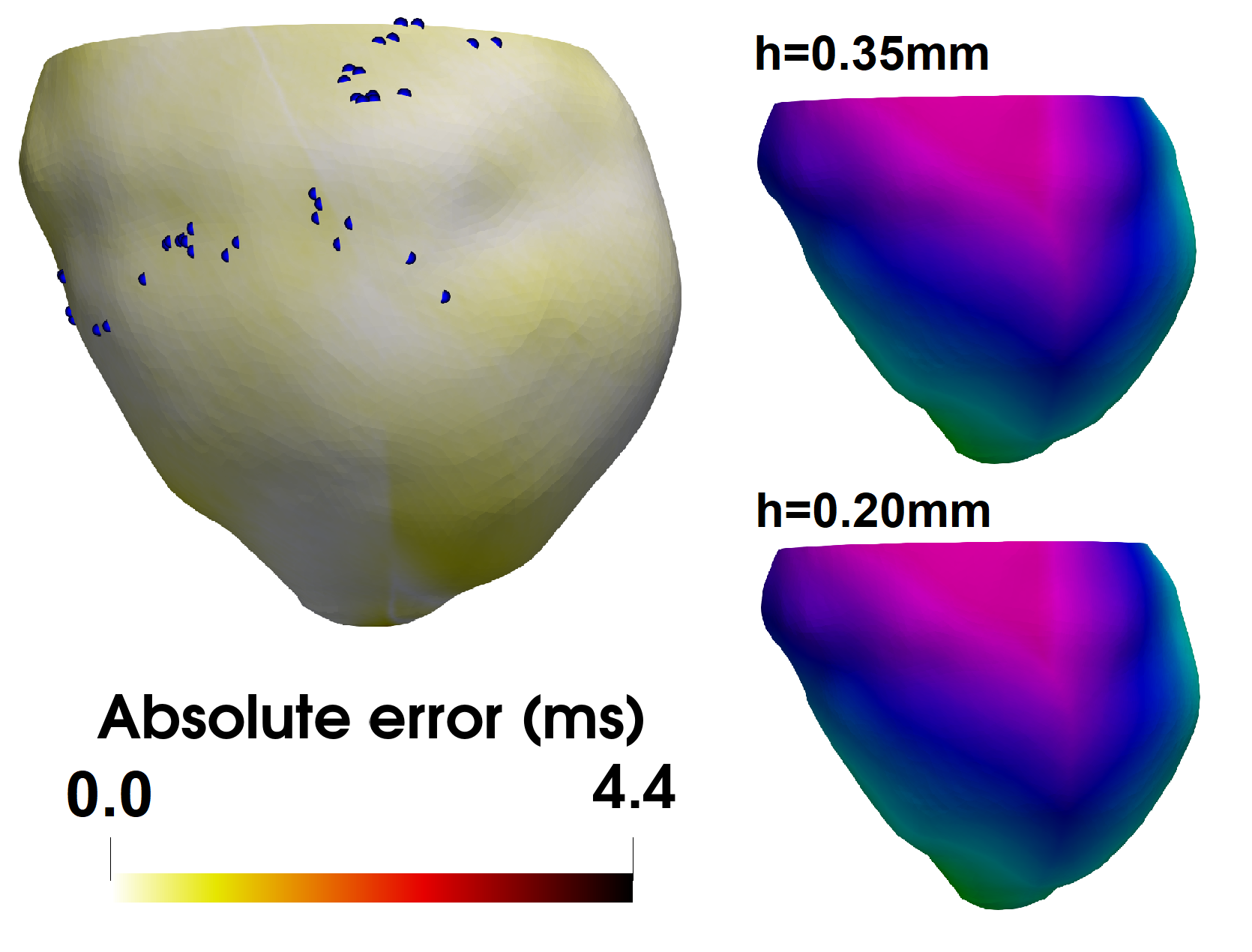}
	\caption{Activation maps obtained for $h\simeq 0.35\,mm$ and $h\simeq 0.20\,mm$, and absolute error
	between them. View of the epicardial veins.} \label{fig:h}
\end{figure}

In this section, we show the numerical results obtained in terms of conductivities estimation 
and corresponding comparison between measured and computed activation times for a standard scenario
obtained by reference numerical and physical parameters, test A in what follows. 
In order to verify that our choice of the mesh size ($h\simeq \,0.35mm$) is adequate, we preliminary ran
a test for P3 also with $h\simeq 0.20\,mm$. In Figure \ref{fig:h} we report the activation maps, the absolute error
and the location of the measured points in the epicardial veins. After a quite small error located at the septum,
the discrepancies between the two cases becomes negligible, in particular at the epicardial veins is almost null.
This allowed to obtain converged results about the conduction velocities in the minimization problem.

The average number of iterations required to solve the optimization problem is 4,
this means that, in average, about 4 monodomain problems need to be solved for each patient
	to estimate his conductivity tensor and each monodomain solution requires about 10 hours of computation. 
	Anisotropy ratio among conductivities was checked a posteriori to fall down in the physiological range.

In Table \ref{table:sigma} we report the estimation of the conductivities obtained for the four patients
by solving the minimization problem \eqref{eq:minimization}
and the mean relative value $\yy e^{I}=\frac{1}{N^{V-I}}\sum_j\frac{e^{I}_j}{\max_i\tau^{V-I}_i}$,
with $e^{I}_j=|\tau^{V-I}_j-\tau^{h-I}_j|$ at the points at the epicardial veins of group I.
Notice that all the values fall in the
physiological range $\Sigma$ reported in Section \ref{sec:invpb}. 
\begin{table}[H]
	\centering
	\begin{tabular}{c|c|c|c|c|}
		\cline{2-5}
		& $\widehat\sigma_f$ & $\widehat\sigma_s$ & $\widehat\sigma_n$ & Mean relative error $e^{I}$ [$\%$] \\
		\hline
		\multicolumn{1}{ |c|   }{P1} & 1.11 & 0.21 & 0.05 & 4.95\\
		\hline
		\multicolumn{1}{ |c|   }{P2} & 1.57 & 0.41 & 0.08 & 5.97 \\
		\hline
		\multicolumn{1}{ |c|   }{P3} & 1.23 & 0.25 & 0.07 & 6.90 \\
		\hline
		\multicolumn{1}{ |c|   }{P4} & 1.39 & 0.30 & 0.07 & 6.12\\
		\hline
	\end{tabular}
	\caption{Values of the optimal conductivity $\widehat{\xx\sigma}$ expressed in $k\Omega^{-1}cm^{-1}$
	and mean relative error in the calibration points of group I. Test A.}\label{table:sigma}
\end{table}
Starting from the conductivity values reported in Table \ref{table:sigma},
	we computed the corresponding planar wave front velocities along the principal axes. We obtained the 
	following average values among the 4 patients: 
	$0.63m/s$, $0.44m/s$ and $0.18m/s$ along the fibers, the transversal and the normal directions, respectively. These values resulted to fall in the physiological ranges, see \cite{Costa_2013,Augustin_2016_bis}.

\begin{figure}[H]
	\centering
	\includegraphics[scale=0.065]{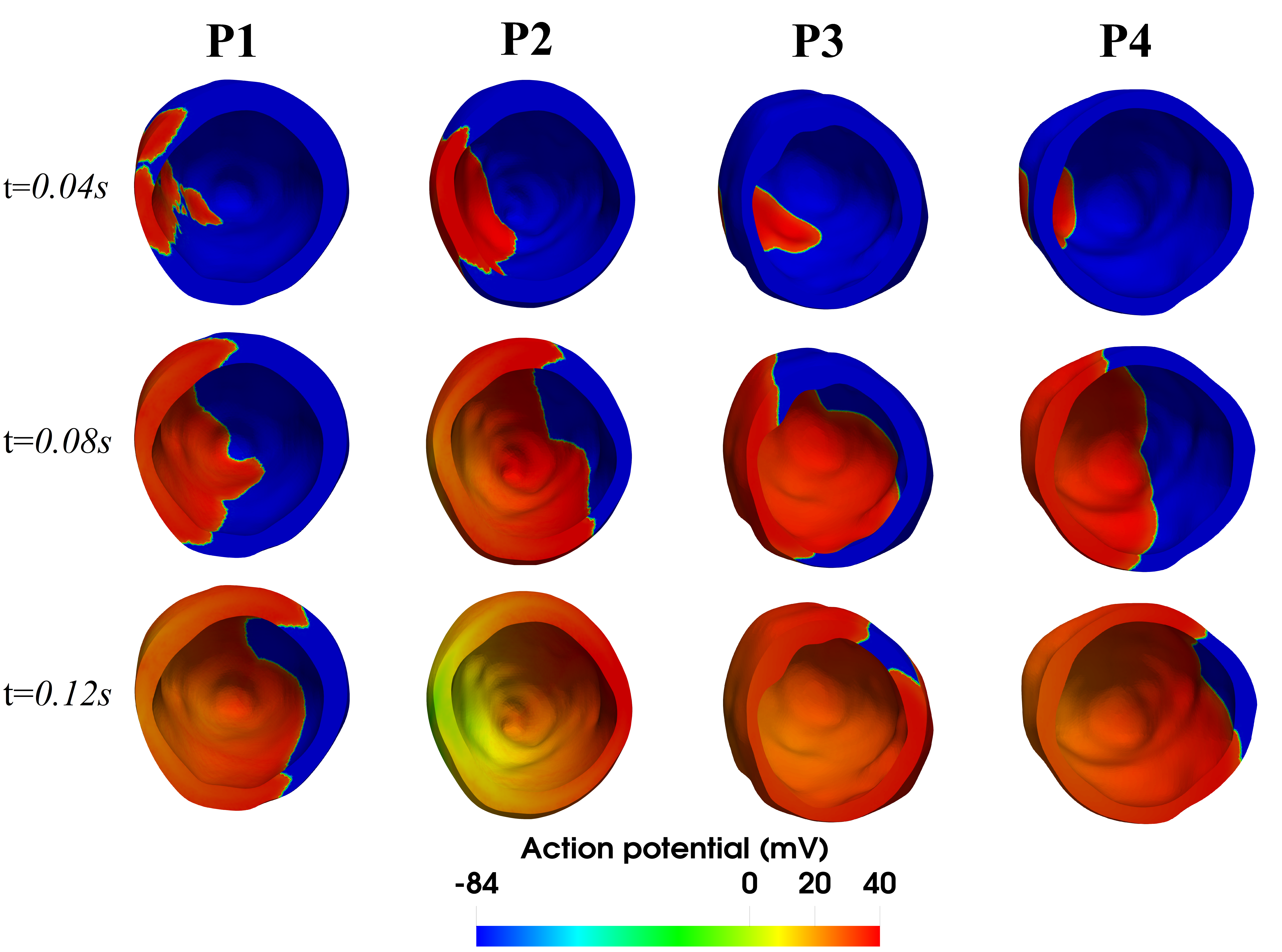}
	\caption{Action potential from the base view of the computational domain. In each row different time instants and in each column different patients are shown. Test A.} \label{fig:5}
\end{figure}

In Figure \ref{fig:5} we show for the four patients the action potential at different instants computed by the numerical simulations with the estimated conductivities reported in Table \ref{table:sigma}.
As highlighted by this figure, the region activated first was the septum, where measurements 
used as input were available.
This is in accordance with the fact that the patients suffer from LBBB, thus the signal does not enter the left 
ventricle through the Purkinje network as in the normal propagation \cite{vigmondc1,bordasg2,vergaral1}, rather through the septum activated by the right ventricle. 

\begin{figure}[H]
	\centering
	\subfloat[\textbf{P1}]{\includegraphics[scale=0.06]{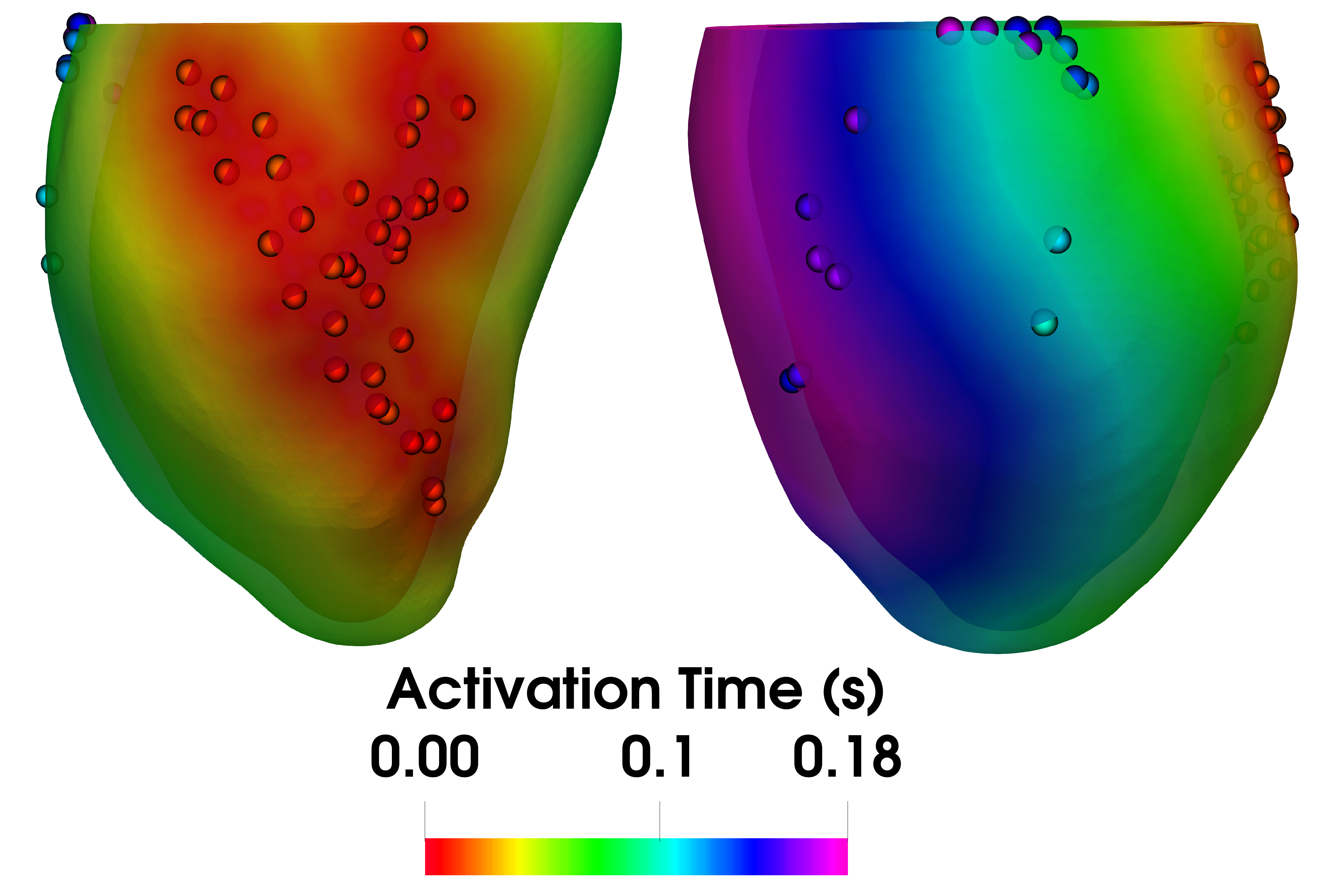}}\qquad\qquad
	\subfloat[\textbf{P2}]{\includegraphics[scale=0.06]{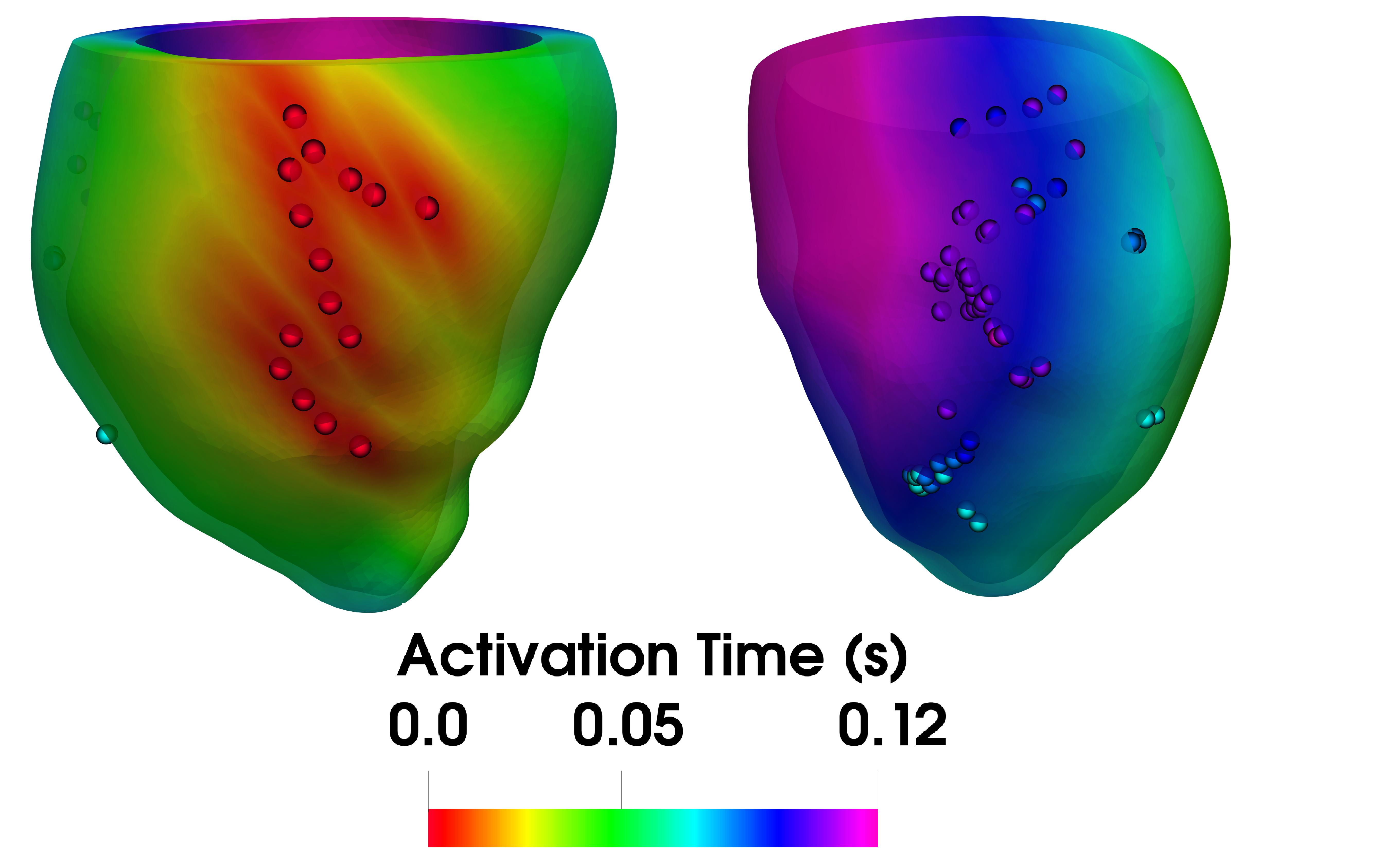}}\\
	\subfloat[\textbf{P3}]{\includegraphics[scale=0.12]{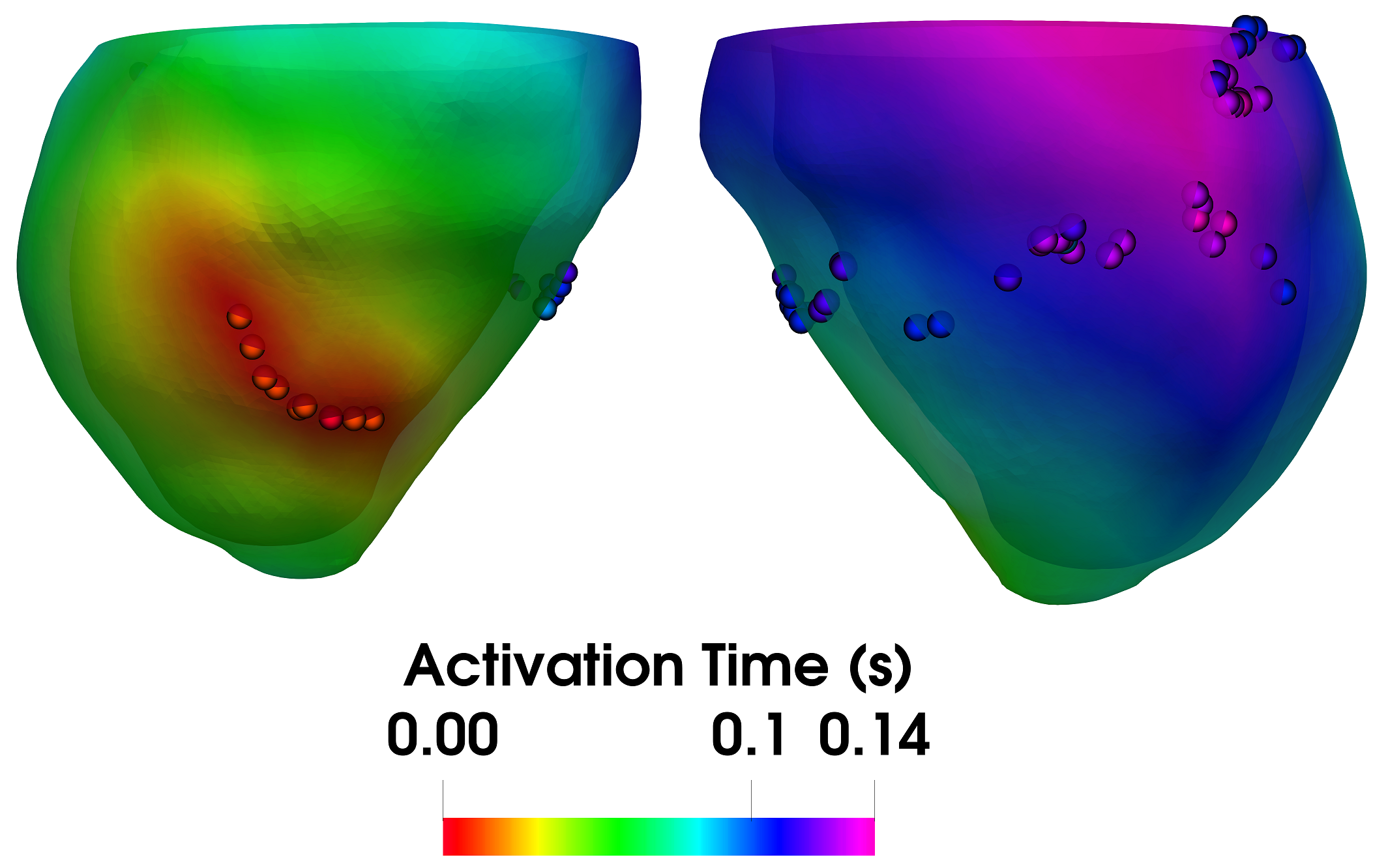}}\qquad
	\subfloat[\textbf{P4}]{\includegraphics[scale=0.08]{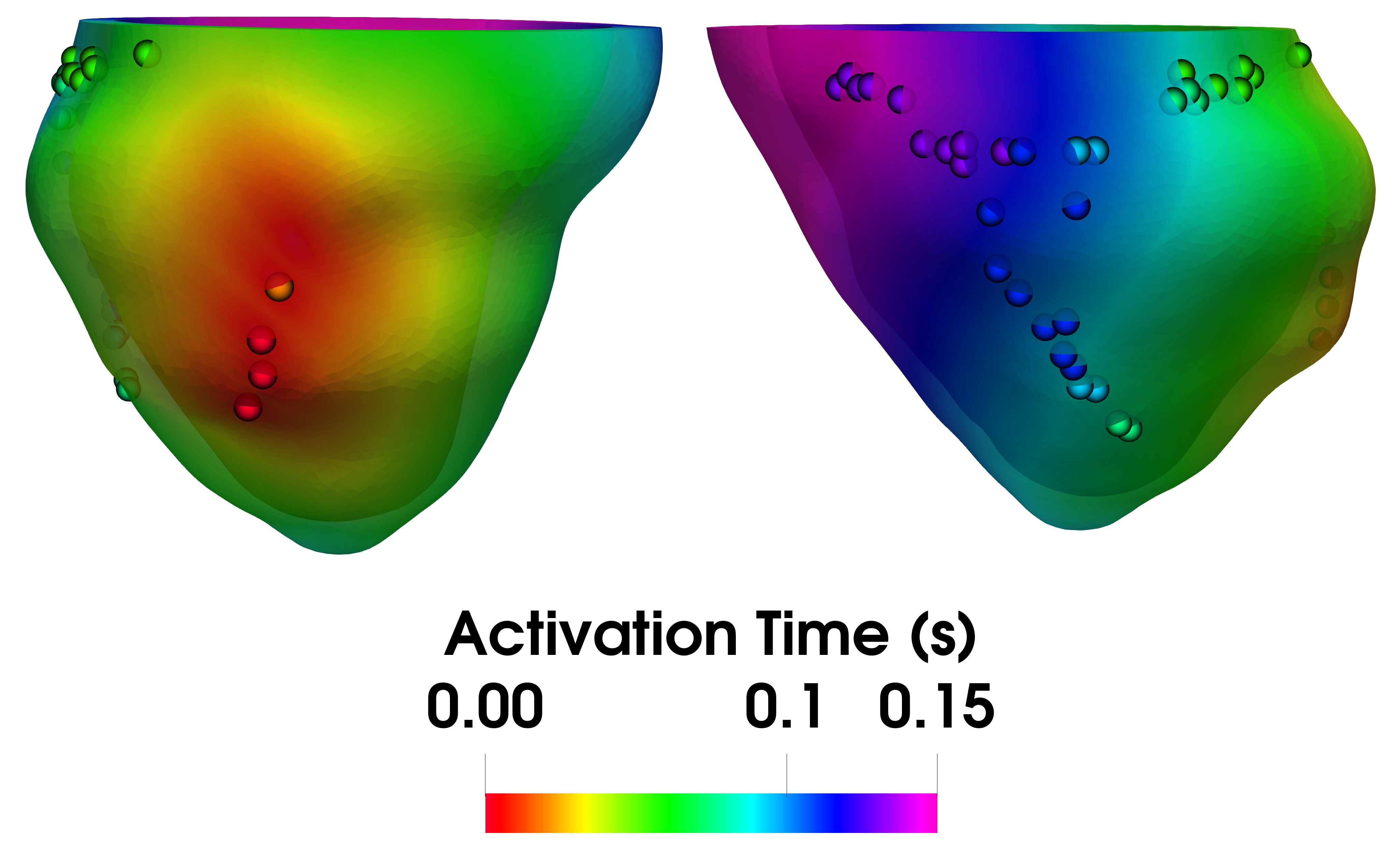}}
	\caption{Computed activation times obtained by numerical solutions (continuous maps) 
		and clinical measurements of activation times (bullets). Test A.} \label{fig:6}
\end{figure}

According to \eqref{eq:tau}, in Figure \ref{fig:6} we report the 
activation times corresponding to the previous numerical results (continuous map) together
with the measurements (bullets).
\captionsetup[subfigure]{labelformat=empty}
From these results, we observe an excellent qualitative agreement between the computed and
measured activation times at the epicardial veins also for group II, that is those
with highest values. 
We remind that the latter measurements were not
used in the minimization problem, thus we have here provided a fair cross validation of the monodomain
problem in the context of sinus rhythm. 

In order to go more in depth in the error analysis, 
in Table \ref{table:errors} we report the values of the errors computed for each patient. 
In particular, we computed the quantities $e^{II}_j=|\tau^{V-II}_j-\tau^{h-II}_j|$
and the corresponding mean relative values $\yy e^{II}=\frac{1}{N^{V-II}}\sum_j\frac{e^{II}_j}{\max_i\tau^{V-II}_i}$ 
and $\yy e^{II,2}=\frac{1}{N^{V-II}}\sum_j\frac{e^{II}_j}{\tau^{V-II}_j}$
together with the corresponding standard deviations over the total number of epicardial veins measurements belonging to group II, see Table \ref{table:measurements}.
\begin{table}[H]
	\centering
	\begin{tabular}{|c||c|c||c|c||c|c|}
		\hline
		&&&&&&\\
		Patient & $e^{II}$ [$\%$] & Std [$\%$] & $e^{II,2}$ [$\%$] & Std,2 [$\%$] &
		Slope $s$ & $R^2$ \\
		\hline
		P1 & 6.17 & 3.28 & 6.67 & 3.11 & 1.0 & 0.65 \\
		\hline
		P2 & 4.12 & 2.19 & 4.41 & 2.07 & 0.96 & 0.72 \\
		\hline			
		P3 & 5.05 & 2.88 & 5.42 & 2.54 & 0.75 & 0.85\\
		\hline
		P4 & 5.42 & 1.95 & 5.75 & 1.92 & 1.11 & 0.89\\
		\hline
	\end{tabular}\caption{Left: Values of the mean relative errors $e^{II}$ and $e^{II,2}$, over the whole set of 
		epicardial veins data belonging to group II, between numerical results and measurements, and corresponding standard deviations. Right: Slope of regression fit $s$ related to correlation plots between numerical
		and measured activation times and coefficient of determination $R^2$. Test A.}\label{table:errors}
\end{table}

\begin{figure}[H]
	\centering
	\includegraphics[width=8cm,height=5.cm]{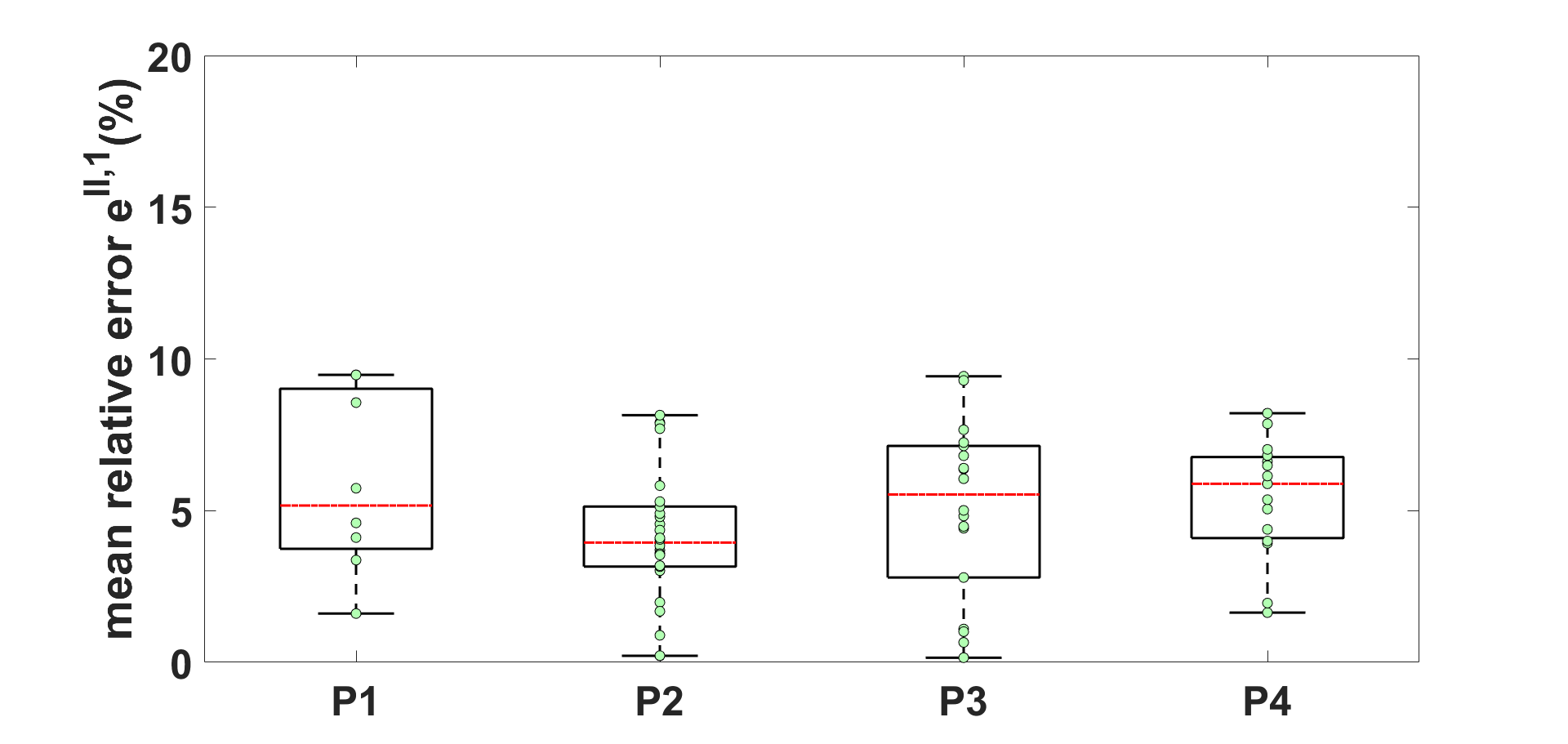}
 \includegraphics[width=8cm,height=5.cm]{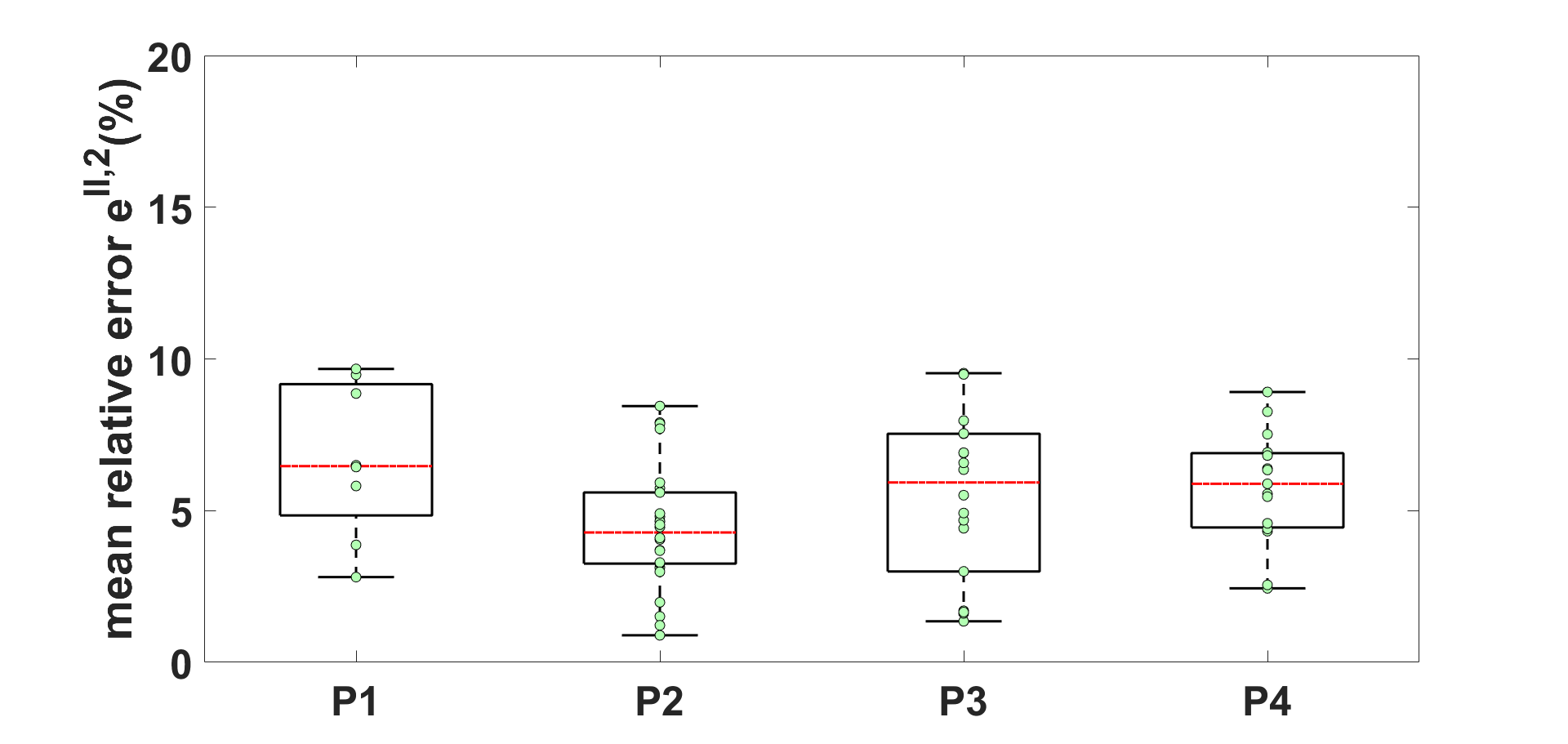}
	\caption{Boxplots of relative errors for the four patients with minimum and maximum values, median,
		and first and third quartiles. Left, the case of error $e^{II}$; Right: the case of error $e^{II,2}$.
		Test A.}
	\label{fig:box}
\end{figure}

For all the four cases, in Figure \ref{fig:box} we report the {\sl boxplots} of the relative errors. This technique is useful to display groups of data through their quartiles. It is based on a five-number summary: the minimum and the maximum values of the dataset (shown by the lower and the upper lines on the whisker), the median (red dashed line), and the first and third quartile (lower and upper bounds of the box). To better show the distribution of the errors, we also reported the values of the single relative error by using green dots. 

From these results, we observe an excellent quantitative agreement between epicardial vein
measurements not used in the parameter estimation (group II) and our computed results. The mean error was 
in any case below $6.2\%$ ($6.7\%$ if $e^{II,2}$ is considered) and the standard deviation confirmed
a low error variability.
Moreover, from the statistical analysis reported in Figure \ref{fig:box},
we can observe that the errors are quite well confined in a small region, the maximum relative error being 
in any case less than 10\%.

\begin{figure}[H]
	\centering
	\includegraphics[scale=0.16]{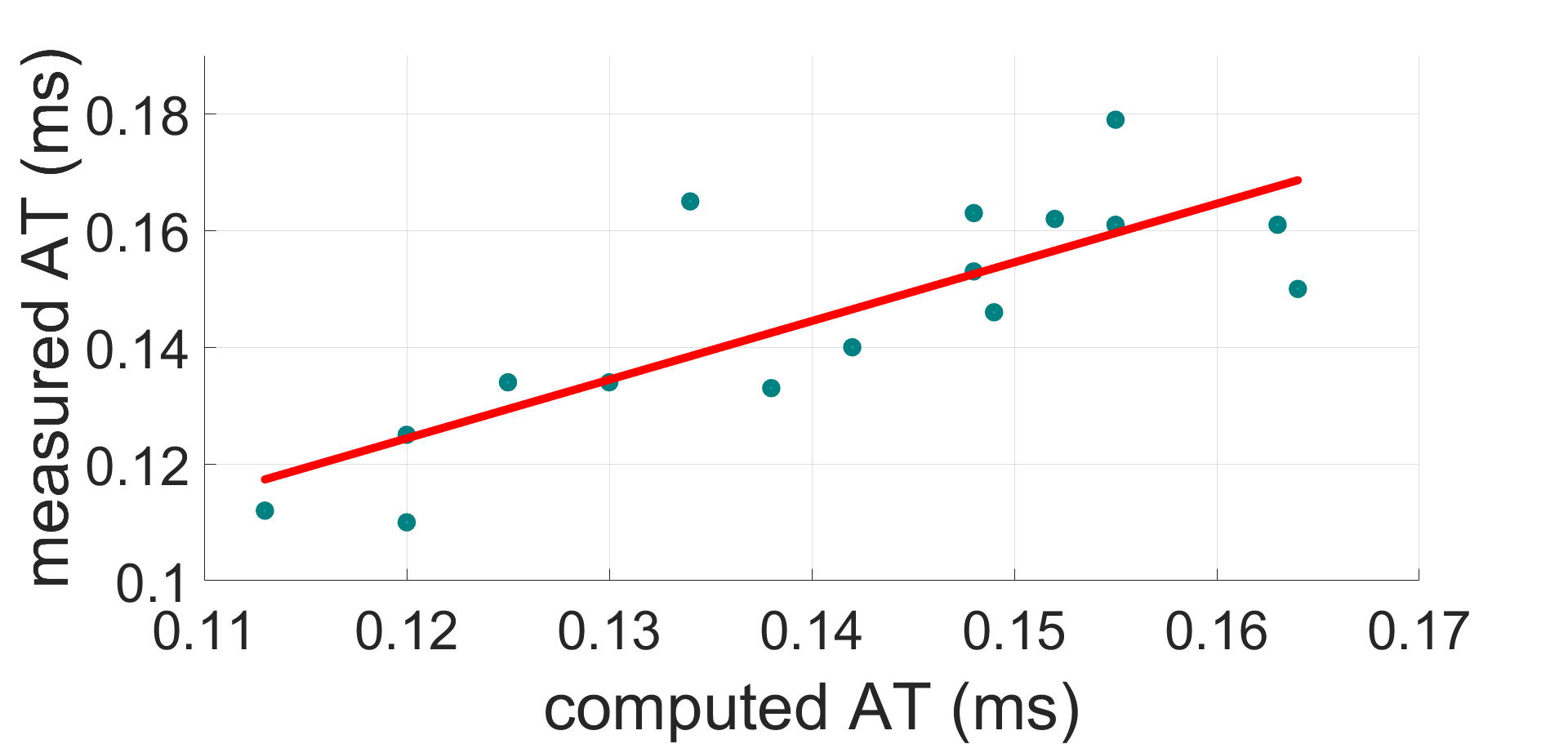}
 \includegraphics[scale=0.16]{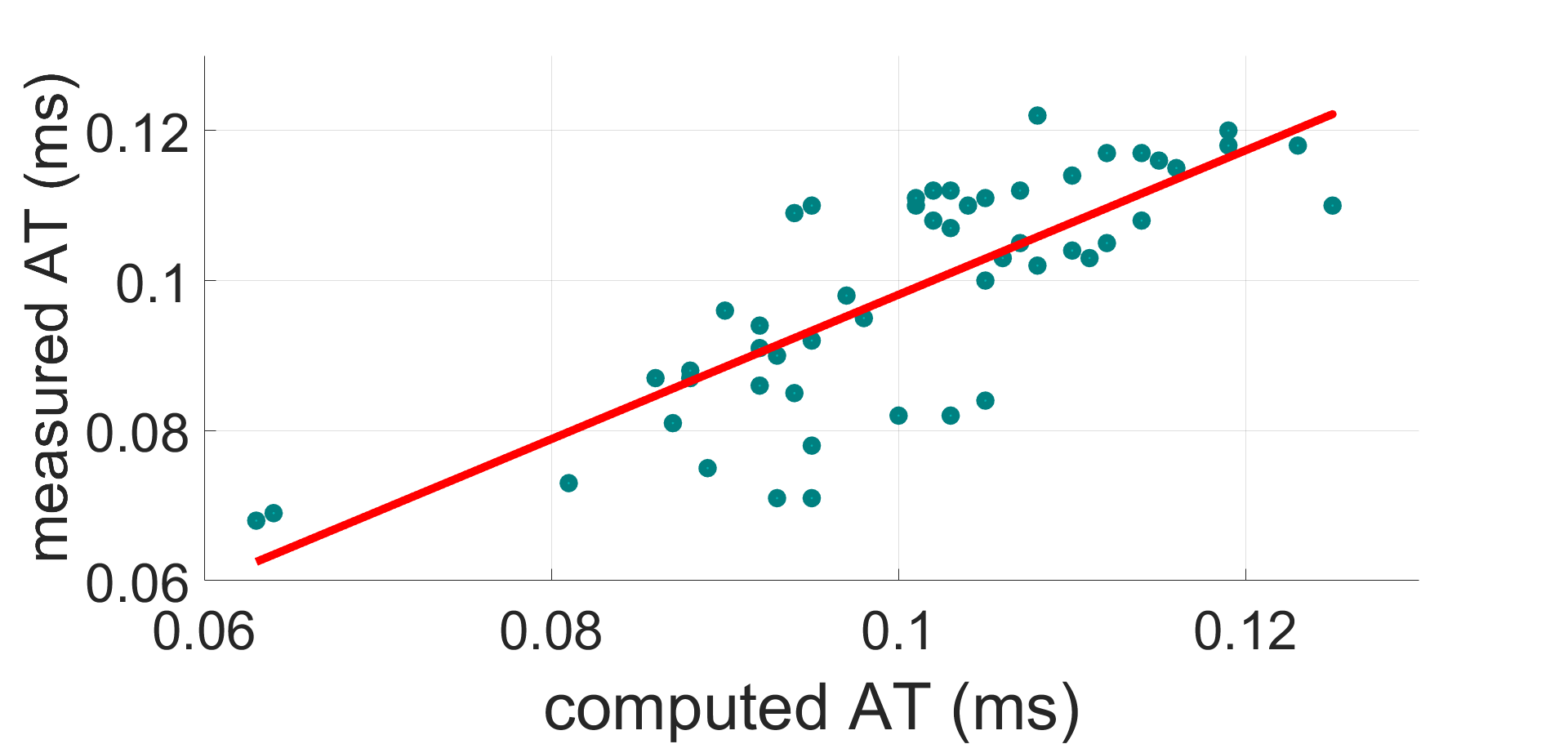}\\
	\includegraphics[scale=0.16]{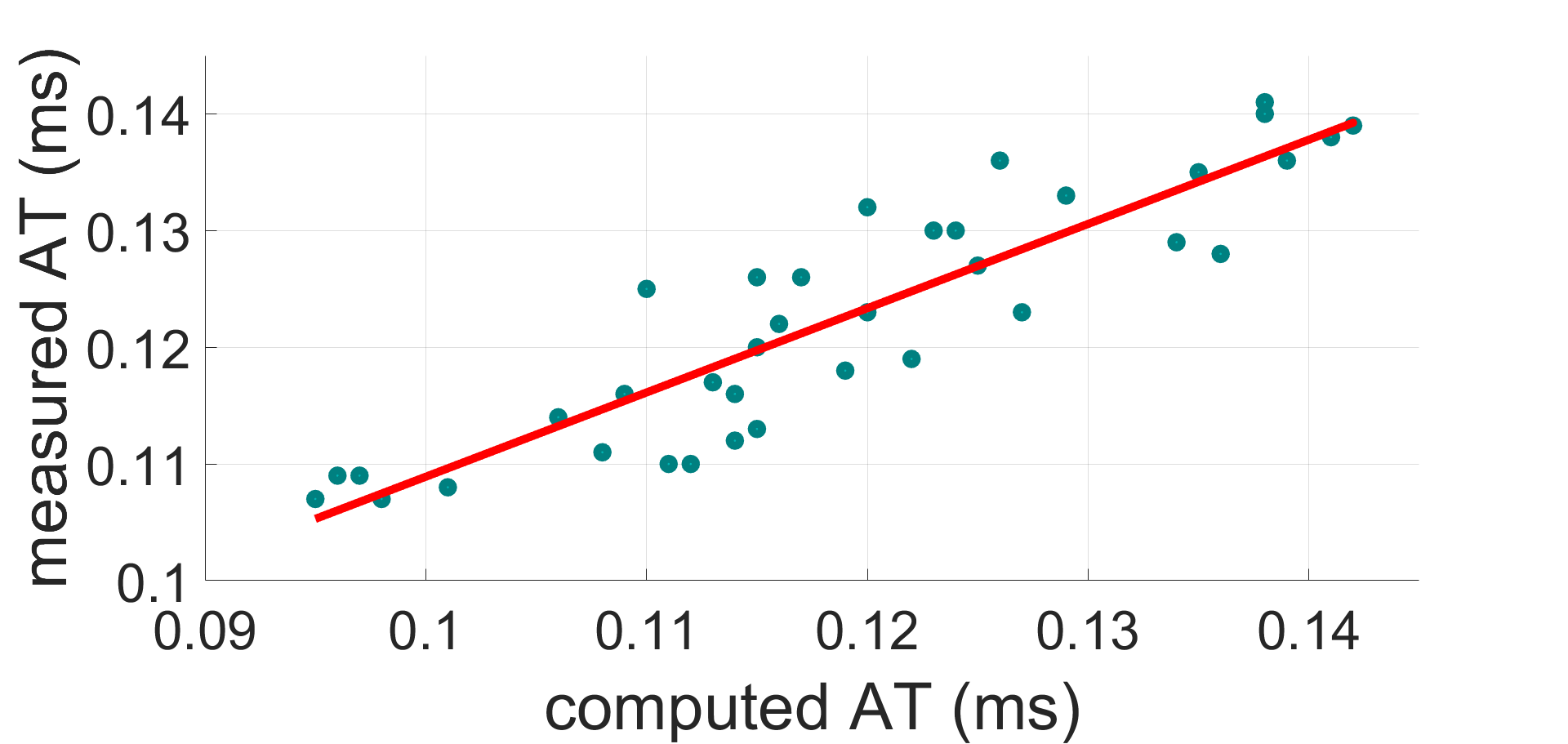}
 \includegraphics[scale=0.16]{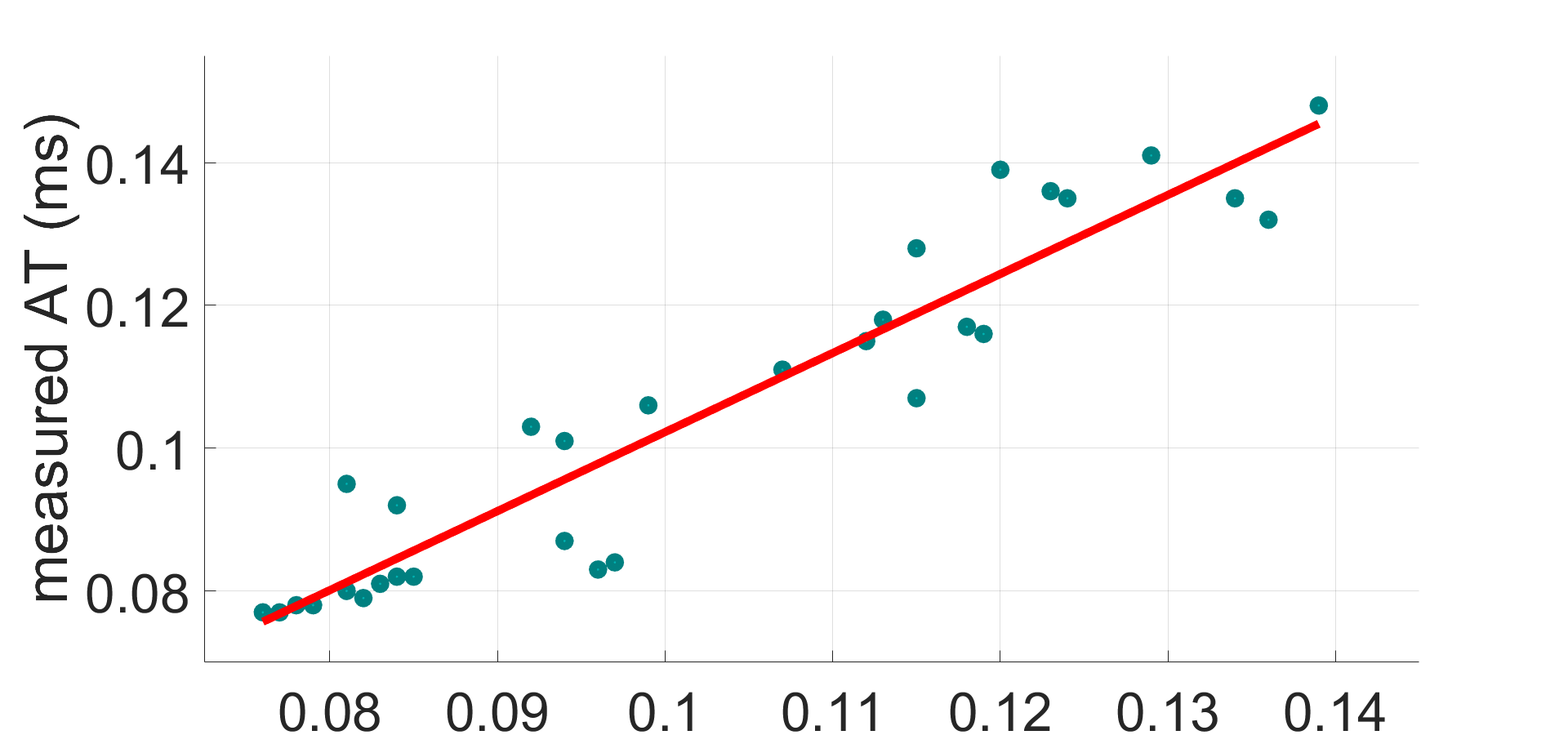}
	\caption{Correlation plots between numerical and measured activation times
	(group I + group II) with regression fitting straight line (in red). Top, left: P1; Top, right: P2; 
	Bottom, left: P3; Bottom, right: P4.}
	\label{fig:corr}
\end{figure}

To complete the analysis, in Figure \ref{fig:corr} we report the correlation plots 
between activation times obtained by numerical simulations and measurements.
The corresponding slope $s$ of the regression fitting straight line and the coefficient of determination $R^2$
are shown in Table \ref{table:errors}. We remember that $R^2\in [0,1]$ provides a measure of how well observed outcomes are replicated by the model \cite{Steel_1960}. $R^2>0.7$ is generally considered strong effect size, whereas
$0.5<R^2<0.7$ is generally considered a moderate effect size \cite{Moore_2013}.	From our results, we observe a slope $s$ close to $1$ for all the patients and a coefficient of determination $R^2$ which features a strong effect
size in 3 of the 4 patients, and a moderate effect size (very close to $0.7$) in 1 case. 

\subsection{Sensitivity of optimization with respect to some parameters}

In order to highlight the sensitivity of our optimization procedure with respect to the choice
of some parameters, we performed some further numerical experiments obtained by considering different scenarios. 
In particular:
\begin{itemize}
\item[-] Test B: Absence of cardiac fibers and different fibers boundary angles;
\item[-] Test C: Different registration between electrical and geometric data;
\item[-] Test D: Different number of points in group I;
\end{itemize}

In order to highlight the importance of including the fibers orientations in the search of
optimal conductivities, we solved the minimization problem \eqref{eq:minimization} for a case without fibers,
where we have only one value $\sigma$ of the conductivity. In Table \ref{table:errors2} we report 
the values of the optimal conductivity and the corresponding relative error with respect to group II for P3. 
From these results
we observe the importance of including a suitable orientation of cardiac fibers to obtain accurate estimation
of the conductivities
\begin{table}[H]
	\centering
	\begin{tabular}{|c|c|c|c||c|c|}
		\hline
		&&&&&\\
		$\widehat\sigma_f$ & $\widehat\sigma_s$ & $\widehat\sigma_n$ & $\widehat\sigma^{no-f}$ & $e^{II}$ [$\%$] &
		$e^{II,no-f}$ [$\%$]\\
		\hline
		1.23 & 0.25 & 0.07 & 1.13 & 5.05 & 7.93\\
		\hline
	\end{tabular}
	\caption{Values of optimal conductivities expressed in $k\Omega^{-1}cm^{-1}$ and  
	mean relative error between numerical results and measurements in presence and absence of cardiac fibers
	(the latter identified by index ''no-f''). P3. Test B.}
	\label{table:errors2}
\end{table}

In order to go deeper in the analysis of the influence of cardiac fibers orientation on the results of the minimization problem, we solved problem \eqref{eq:minimization} for other two boundary values for the fibers
in the Laplace-Dirichlet rule-based algorithm, that is $\pm 45^\circ$
and $\pm 75^\circ$, the values of the sheets angles being the same of above. 
In Table \ref{table:errors3} we reported the optimized conductivities and mean relative error for P3. 
From these results we observe that the value $60^\circ$ produces the smallest errors among the three choices.
\begin{table}[H]
	\centering
	\begin{tabular}{|c||c|c|c||c|}
		\hline
		&&&&\\
Fibers angle &		$\widehat\sigma_f$ & $\widehat\sigma_s$ & $\widehat\sigma_n$ & $e^{II}$ [$\%$] \\
		\hline
$\pm 45^\circ$ & 1.12 & 0.22 & 0.07 & 7.67\\
$\pm 60^\circ$ & 1.23 & 0.25 & 0.07 & 5.05\\
$\pm 75^\circ$ & 1.36 & 0.27 & 0.07 & 7.55\\
		\hline
	\end{tabular}
	\caption{Values of optimal conductivities expressed in $k\Omega^{-1}cm^{-1}$ and  
	mean relative error between numerical results and measurements for different fibers boundary angles. P3.
	Test B.}
	\label{table:errors3}
\end{table}

\begin{figure}[H]
	\centering
	\includegraphics[scale=0.15]{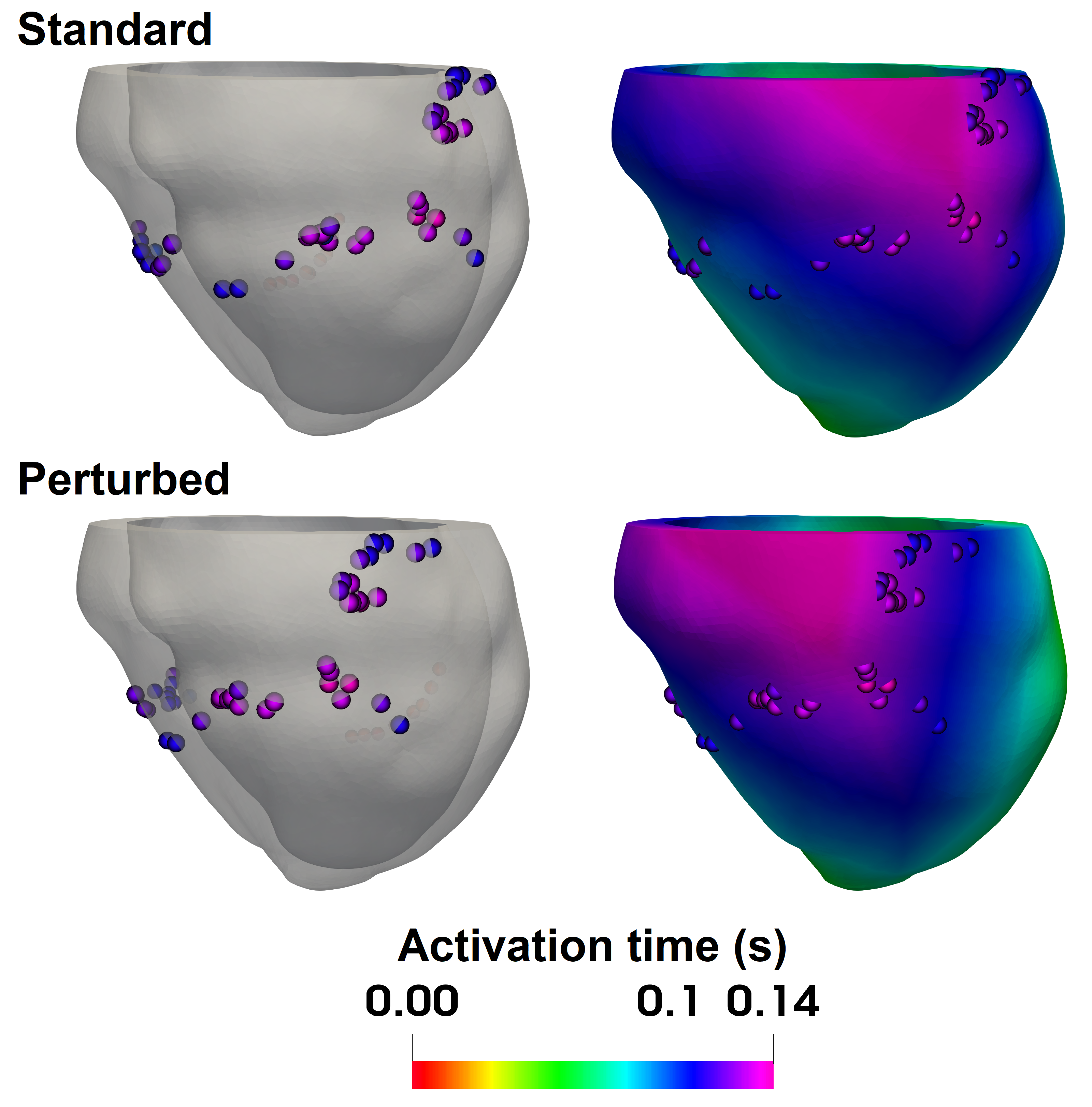}
	\caption{Different merging between electrical and geometric data. Left: Location of the measurements in the epicardial vein; Right: activation maps from
	numerical simulations. Top: standard scenario; Bottom: perturbed scenario obtained by changing the location
	of the reference point during the merging of electrical and geometric data. P3. Test C.}
	\label{fig:registr}
\end{figure}

As observed in Sect. 2.2.2, the merging of electrical and geometric data 
was based on the a priori selection of some reference points. This choice was guided by the anatomy of the patient
(coronary sinus and septum) but it could be affected at a certain extent by an arbitrariness 
having an impact on the results. To investigate this point, we considered for P3 another configuration
({\sl perturbed scenario}),
where the reference points were located as much as possible far from the original ones, although  
satisfying anatomical requirements, see Figure \ref{fig:registr}, left.

As observed by the activation maps reported in Figure \ref{fig:registr}, right,
there is a good agreement between measured and computed data also in the perturbed scenario.
In fact, the different location of the septal reference points induces a rigid translation of all the other
measures, so that it is not surprising that the qualitative propagation map is still well
fitting the data. 
This is confirmed by the errors for data groups I and II reported in Table \ref{table:errors4},
which show almost identical values, suggesting the robustness of our analysis with respect
to the arbitrariness of the data merging procedure.
\begin{table}[H]
	\centering
	\begin{tabular}{|c||c|c|}
		\hline
		&&\\
Scenario &		$e^{I}$ [$\%$] & $e^{II}$ [$\%$]\\
		\hline
		\hline
 standard & 6.90 & 5.05 \\
 geometric perturbed & 6.84 & 5.10  \\
		\hline
		$N^V_I=15$ & 6.59 & 5.11  \\ 
		$N^V_I=12$ & 6.32 & 5.15  \\
		\hline
	\end{tabular}
	\caption{Relative errors in correspondence of the two groups of data (I and II)
	for standard scenario, geometric perturbed scenario (Test C) and different number of points in group I (Test D). P3.}
	\label{table:errors4}
\end{table}

We want now to assess the importance of the number of points of group I used
to calibrate the conduction velocities on the accuracy of the optimization procedure. 
To this aim, we considered for P3 2 further scenarios obtained by using 15 and 12 points
(19 were used in the standard scenario). In Table \ref{table:errors4} we reported the corresponding errors,
obtained after the same optimization procedure that led to the same values of conductivities of the standard scenario.
These errors highlight the robustness of the method with respect to the number of mapped points.  

%%%%%%%%%%%%%%%%%%%%%%%%%%%%%%%%%%%%%%%%%%%%%%%%%%%%%%%%%%%%%%%%%%%%%%%%%%%%%%%%%%%%%%%%%%%%%%%%%%%%%

%%%%%%%%%%%%%%                             SECTION 4                         %%%%%%%%%%%%%%%%%%%%%%%%
\section{Discussion}

\subsection{General discussion of the results}
In this paper we have performed an important step towards the validation of the monodomain model 
in the context of a normal ventricular sinus rhythm (i.e. without chaotic patterns) by using in-vivo data of activation times acquired in four patients. At the best of authors knowledge, this is the first validation
attempt of the monodomain model where measured ventricular activation maps in humans are used to estimate conduction parameters and assess the errors in a cross-validation test.   

Specifically, for each of the four cases considered, we have used patient-specific electrical 
data acquired at the septum to provide an input which is suitable for LBBB patients. 
We used	a set of measures 
at the epicardial veins to calibrate the three conductivities and another one 
at the same location to quantify the error of the numerical results. The mean error found in this cross-validation test
was in any case less than 6.2\%, see also Figures \ref{fig:6}, \ref{fig:box} and Table \ref{table:errors},
and was comparable with that obtained for group I (see Table \ref{table:sigma}).This proves the robustness of the minimization problem for the data we had at disposal.  
From Table \ref{table:measurements} we also notice that for P2 and, especially, for P3 and P4
our accurate results were found by using few input data (septal measurements). This means that it is enough to have
at disposal very few measurements of activation times where the signal starts in order to well predict the
activation in the epicardial veins. 

The validation against epicardial measures is particularly significant. Indeed,
the front propagating from the endocardium activates all the myocardium before reaching the epicardium,
thus assessing the accuracy of the results at the epicardium means accounting for the reliability of the 
propagation in all the myocardium. Among the studies which incorporates epicardial measures in computational models, only \cite{Lee_2019}, together with the present study, considered measures at the epicardial veins. 
Some differences of the present work with respect to \cite{Lee_2019} 
are the use of monodomain model instead of the reaction-Eikonal one, the simulation of a sinus rhythm 
instead of a stimulated scenario (at the right ventricle apex), the use of patient specific input data at the septum,
and the calibration against epicardial activation times instead of ECG measures.
 
Our modeling choice was oriented towards the monodomain model, at the expense of an increased computational time
with respect, e.g., to the Eikonal models. 
	This was motivated by the fact that the monodomain model should be more accurate.
	Indeed, although the Eikonal equation was seen to be accurate in predicting the activation time, it is a reduced model of the monodomain one \cite{Colli_1990} and thus introduces an error
	in the computation of the activation times with respect to monodomain. For example, in \cite{Wallman_2012},
	such discrepancy has been quantified to be about 10\%.
	Moreover, the monodomain model allowed us to compute also the transmembrane potential
	and thus to provide supplementary clinical information. Additionally, the monodomain model will be fundamental 
	for electro-mechanical simulations to effectively assess the functioning of CRT.
An alternative choice could be provided by the reaction-Eikonal model \cite{Neic_2017,Lee_2019},
which at some extent merges Eikonal and monodomain models. 
This method has a good accuracy (discrepancies with respect to the bidomain model in any case less than 20\%,
with about 33\% of the points affected by an error smaller than 2.5\%) and very good computational times, 
which are intermediate between Eikonal and monodomain. Moreover, unlike Eikonal this model is able to provide the 
transmembrane potential. It will however be interesting for further studies to assess if any differences 
between reaction-Eikonal and monodomain model would be significantly out-weighted by other errors in the simulation pipeline 
(registration errors, inaccurate representation of fibre architecture, for example). 

Results in Table \ref{table:errors2} and \ref{table:errors3} showed the importance of including a suitable cardiac fibers orientation
when looking for optimal values of the conductivities to accurately match the clinical measurements. 
In particular, although $45^\circ$ and $75^\circ$ are acceptable boundary values for fibers angle
falling in the physiological range \cite{Lombaert_2012}, the value $60^\circ$ produced the smallest errors,
confirming the validity of this choice as proven by the hystological studies in \cite{Toussaint_2013}.
Moreover, the absence of any fibers, due to the simplified model (only one parameter instead of three to estimate)
led to error about 60\% greater than in the case with a suitable cardiac fibers field.

Some of the results presented are, at the best of authors knowledge, completely new in the field
of computational electro-physiology. First, we highlight that this is the first cross-validation analysis against clinical measures of ventricular activation times. Second, for the first time we used measures at the epicardial veins
to validate the monodomain equation. This is particularly relevant since the signal usually 
originates at the endocardial level. Finally, we included some results about the sensitivity of the solution
of the minimization problem with respect to the cardiac fibers orientation.

\subsection{Possible clinical implications}
The previous findings could be of some interest for clinical applications.
For example, in view of CRT cardiologists often use the point 
with the latest activation time (LAT) at the epicardial veins to locate the left electrode
\cite{Zanon_2014,Liang_2015}.
Then, thanks to our accurate numerical method, this information can be in principle provided
without a complete mapping of the epicardial veins, shortening the invasive procedure 
based on the insertion of catheters (up to a couple of hours) and
reduce the exposition of the patient to radiation.  
For example, only the mapping of the coronary sinus could be in principle provided
to obtain the data needed for the calibration. This mapping is in general quite ''simple'' to perform
since the coronary sinus is easily detectable and in the most proximal region.

As observed, the choice of using the earliest (proximal) and latest (distal) activated points for the calibration and validation groups, respectively, allowed us to avoid overlaps and obtain more significant results.
Moreover, this could be of particular clinical interest, 
since one could acquire only the earliest (most accessible) activated points, e.g. those in the coronary sinus, 
and computationally estimate (without the need of mapping) the deepest ones, where LAT occurs.
In this way there is no need to push the catheter deeper in the veins, where the activation times could be completely predicted 
by our computational method, yielding a potential benefit for the patient.
This has been emphasized from the analysis of Test C (see Table \ref{table:errors4})
where we stressed the need of few points (for example the more proximal ones on the coronary sinus) 
to be acquired by the mapping procedure in order to have an accurate calibration.  
Moreover, all the epicardial surface could be virtually mapped.
This would allow to find the global LAT (i.e. the latest among all the epicardial points, not only those in the veins) which could be useful for surgical implantation of the left
electrode \cite{Navia_2005,Marini_2019}.

The effective applications of the monodomain model for such purposes
is however nowadays a little bit problematic since the implantation of CRT and the mapping procedure
happen during the same day
and computational costs of monodomain are still of several hours, thus requiring high performance computing resources.
However, this strategy could be effective in few years owing to the development of more performing computing resources, 
new preconditioners, and implementation of faster optimization strategy.

Another important possible clinical implication made possible by the use of a validated monodomain model
is the electrical-mechanical (EM) virtual study of CRT where the left 
electrode is placed on the computed LAT. Nowadays, this application would be limited to case studies due to the high computational time of EM simulations. However, this could provide
useful information about the effective validity of LAT as stimulation point able to provide a restored heart function. 

\subsection{Limitations}
Some limitations affected the present study. 

The first one consisted, as observed, 
in the absence of any model for the Purkinje system, both in terms of reconstruction of a suitable network and
by surrogating its effect increasing the LV subendocardium conduction \cite{Draper_1959},
as done first in \cite{Cardone_2016} and then also, e.g., in \cite{Lee_2019}.

Another limitation of the current work consists in the comparison of numerical results
only against sparse activation points. A more complete validation could in future be obtained by using full endocardial or epicardial maps. This however will require specific acquisition protocols since full maps are usually outside 
the standard clinical practice. Also, when available, the comparison with repolarization maps and not only 
with activation times will be mandatory to complete the validation of monodomain model.

We also mention that here we focused only on data of sinus rhythm, 
thus we did not validate restitution properties of the model and pathological scenarios. 
Cross-validations against measures obtained under suitable stimulus protocols and for arrhythmic cases
are mandatory to complete our work.

Another limitation is the use of the left ventricle solely, instead of a biventricular model
as done in \cite{Lee_2019}. The inclusion of the right ventricle will be mandatory in electro-mechanical studies
for CRT applications.  
	We notice however that, since we are interested here in predicting the results in the left ventricle epicardium, the absence of the right ventricle geometry should not affect the results. This was motivated by the fact that we used patient-specific electrical data at the septum as input for the left ventricle activation.

We finally mention the fact that we did not address the issue of quantifying the geometric projection error obtained when merging the electrical data from Ensite Precision system with the epicardial surface of the ventricle, see Figure \ref{fig:2}. Currently, we are investigating how much this error can affect our error analysis.

\section*{Acknowledgements}
This project has received funding from the European Research Council (ERC)
under the European Union’s Horizon 2020 research and innovation programme (grant
agreement No 740132, iHEART - An Integrated Heart Model for the simulation of the
cardiac function, P.I. Prof. A. Quarteroni).
CV has been partially supported also by the H2020-MSCA-ITN-2017, EU project 765374
"ROMSOC - Reduced Order Modelling, Simulation and Optimization of Coupled systems".
AQ and CV have been partially supported also by
the Italian research project MIUR PRIN17 2017AXL54F "Modeling the heart across the scales:
from cardiac cells to the whole organ".

%%%%%%%%%%%%%%%%%%%%%%%%%%%%%%%%%%%%%%%%%%%%%%%%%%%%%%%%%%%%%%%%%%%%%%%%%%%%%%%%%%%%%%%%%%%%%%%%%%%%%

%%%%%%%%%%%%%%                           BIBLIOGRAPHY                        %%%%%%%%%%%%%%%%%%%%%%%%
\bibliographystyle{plain}
\bibliography{biblio}
%%%%%%%%%%%%%%%%%%%%%%%%%%%%%%%%%%%%%%%%%%%%%%%%%%%%%%%%%%%%%%%%%%%%%%%%%%%%%%%%%%%%%%%%%%%%%%%%%%%%%

\end{document}